\documentclass[twoside,11pt]{amsart}
\usepackage[active]{srcltx} 
\usepackage[all]{xy}
\CompileMatrices
\usepackage{amsfonts}
\usepackage{amssymb}
\usepackage{amsthm}
\usepackage{amsmath}
\usepackage{ifthen}
\usepackage{url}
 \newtheorem{thm}{Theorem}[section]
 \newtheorem{prop}[thm]{Proposition}
 
 \newtheorem{lem}[thm]{Lemma}
  
 \newtheorem{conj}[thm]{Conjecture}

\theoremstyle{definition}

\newtheorem{ex}[thm]{Example}

\theoremstyle{remark}
\newtheorem{rem}[thm]{Remark}

\font\russ=wncyr10 \font\Russ=wncyr10   scaled\magstep 1
\def\sha{\hbox{\russ\char88}}
\def\Sha{\hbox{\Russ\char88}}

\renewcommand{\L}{\ifmmode {\mathcal{L}}\else$\mathcal{L}$\ \fi}

\renewcommand{\d}{\mathbf{d}}
\renewcommand{\u}{\mathbf{1}}
\newcommand{\id}{\mathrm{id}}
\newcommand{\bbC}{\ifmmode {\mathbb{C}}\else$\mathbb{C}$\ \fi}
\newcommand{\bbR}{\ifmmode {\mathbb{R}}\else$\mathbb{R}$\ \fi}
\newcommand{\R}{\mbox{$\Bbb R$}}
\renewcommand{\r}{\mathrm{R}\Gamma}

\newcommand{\be}{\begin{equation}}
\newcommand{\ee}{\end{equation}}
\newcommand{\lnr}{ \widehat{L^{nr}}}
\newcommand{\qnr}{ \widehat{\mathbb{Q}_p^{nr}}}
\newcommand{\zpnr}{ \widehat{\mathbb{Z}_p^{nr}}}
\newcommand{\fpbar}{\ifmmode {\overline{\mathbb{F}_p}}\else$\mathbb{F}_p$\ \fi}

\newcommand{\fp}{\ifmmode {\mathbb{F}_p}\else$\mathbb{F}_p$\ \fi}
\newcommand{\zp}{\ifmmode {\mathbb{Z}_p}\else$\mathbb{Z}_p$\ \fi}

\newcommand{\z}{\mathbb{Z}}

\newcommand{\zpMod}{\ifmmode\mbox{$\zp$-Mod}\else$\zp$-Mod \fi}
\newcommand{\Mod}{\ifmmode\mbox{$\Lambda$-Mod}\else$\Lambda$-Mod \fi}
\renewcommand{\mod}{\ifmmode\mbox{$\Lambda$-mod}\else$\Lambda$-mod
\fi}
\newcommand{\La}{\ifmmode\Lambda\else$\Lambda$\fi}
\newcommand{\rhom}{\mbox{$\mathbf{R}\mbox{Hom}$}}
\newcommand{\Hom}{{\mathrm{Hom}}}

\newcommand{\rk}{{\mathrm{rk}}}

\renewcommand{\H}{\mathrm{H}}

\newcommand{\M}{\ifmmode {\frak M}\else${\frak M}$ \fi}
\newcommand{\m}{\ifmmode {\frak m}\else$\frak m$ \fi}
\newcommand{\mh}{\ifmmode {\frak m}(H)\else${\frak m}(H)$ \fi}
\newcommand{\p}{\ifmmode {\frak p}\else${\frak p}$\ \fi}
\renewcommand{\P}{\ifmmode {\frak P}\else${\frak P}$\ \fi}
\newcommand{\e}{\ifmmode {\mathcal{E}}\else$\mathcal{E}$ \fi}

\newcommand{\h}{{\mathcal{ H}}}
\newcommand{\C}{\mathcal{C}}
\newcommand{\B}{\mathcal{B}}
\newcommand{\T}{\mathbb{ T}}
\renewcommand{\O}{\mathcal{ O}}
\newcommand{\G}{\ifmmode {\mathcal{G}}\else${\mathcal{G}}$\ \fi}
\newcommand{\A}{\ifmmode {\mathcal{A}}\else${\mathcal{ A}}$\ \fi}

\renewcommand{\projlim}[1] {{\lim\limits_{\stackrel{\displaystyle
\longleftarrow}{#1}}}}

\newcommand{\kl}{[\![}
\newcommand{\kr}{]\!]}
\newcommand{\Qp}{\ifmmode {{\Bbb Q}_p}\else${\Bbb Q}_p$\ \fi}
\newcommand{\qp}{\ifmmode {{\Bbb Q}_p}\else${\Bbb Q}_p$\ \fi}
\newcommand{\ql}{\ifmmode {{\Bbb Q}_l}\else${\Bbb Q}_l$\ \fi}
\newcommand{\Q}{\ifmmode {\Bbb Q}\else${\Bbb Q}$\ \fi}
\newcommand{\q}{\ifmmode {\Bbb Q}\else${\Bbb Q}$\ \fi}

\newcommand{\coker}{\mathrm{coker}}

\def\sectionnam{\@empty}
\def\subsectionnam{\@empty}

\begin{document}

\title[BSD, ETNC and MC]{From the Birch \& Swinnerton-Dyer Conjecture  over the Equivariant
Tamagawa Number Conjecture to  non-commutative Iwasawa theory -  {a survey}\\
\
\\ {\small\rm after Burns/Flach, Fukaya/Kato, Huber/Kings, Coates, Sujatha ...}\\ \ \\}%

\author{Otmar Venjakob}%
\address{Universit\"{a}t Heidelberg\\ Mathematisches Institut\\
Im Neuenheimer Feld 288\\ 69120 Heidelberg, Germany.} \email{otmar@mathi.uni-heidelberg.de}
\urladdr{http://www.mathi.uni-heidelberg.de/\textasciitilde otmar/}
\thanks{This survey was written during a stay at   Centro de Investigacion en Matematicas (CIMAT), Mexico, as a Heisenberg fellow of the Deutsche Forschungsgemeinschaft (DFG) and I
want to thank these institutions  for their hospitality and
financial support.}



\date{\today}%
\maketitle
\thispagestyle{empty}

This paper  aims to give a survey on Fukaya and Kato's article
\cite{fukaya-kato}  which establishes the relation between the
Equivariant Tamagawa Number Conjecture (ETNC) by Burns and Flach
\cite{bf} and the noncommutative Iwasawa Main Conjecture (MC)
(with $p$-adic $L$-function) as formulated by Coates, Fukaya,
Kato, Sujatha and the author \cite{cfksv}. Moreover, we
compare their approach with that of Huber and Kings \cite{hu-ki}
who formulate an Iwasawa Main Conjecture   (without $p$-adic
$L$-functions). We do not discuss these conjectures in full
generality here, in fact we are mainly interested in the case of
an abelian variety defined over $\Q.$ Nevertheless we   formulate
the conjectures for general motives over $\Q$ as far as possible.
We follow closely the approach of Fukaya and  Kato but our
notation is sometimes inspired by \cite{bf,hu-ki}. In particular,
this article does not contain any new result, but hopefully serves
as introduction to the original articles. See \cite{ven-stockholm}
for a more down to earth introduction to the $GL_2$ Main
Conjecture for an elliptic curve without complex multiplication.
There we had pointed out that the Iwasawa main conjecture for an
elliptic curve is morally the same as the (refined) Birch and
Swinnerton Dyer (BSD) Conjecture for a whole tower of number
fields. The work of Fukaya and Kato makes this statement precise
as we are going to explain in these notes. For the convenience of
the reader we have given some of the proofs here which had been
left as an exercise in \cite{fukaya-kato} whenever we had the
feeling that the
presentation of the material becomes more transparent thereby.\\
Since the whole paper bears an expository style we omit a lengthy
introduction and just state briefly the content of the different
sections:

In section 1 we recall the fundamental formalism of
(non-commutative) determinants which were introduced first  by
Burns and Flach to formulate equivariant versions of the TNC. In
section 2 we briefly discuss the setting of (realisations of)
motives as they are used to formulate the conjectures concerning
their $L$-functions, which are defined in section 3. There, also
the absolute version of the TNC is discussed, which  predicts the
order of vanishing of the $L$-function at $s=0,$ the rationality
and finally the precise value of the leading coefficient  at $s=0$
up to the period and regulator. In subsection 3.1 we sketch how
one retrieves the BSD conjecture in its classical formulation if
one applies the TNC to the motive $h^1(A)(1)$ of an abelian
variety $A/\Q.$ Though well known to the experts this is not very
explicit in the literature. In section 4 we consider a $p$-adic
Lie extension of $\Q$ with Galois group $G.$ In this context the
TNC is extended to an equivariant version using the absolute
version for all twists of the motive by certain representations of
$G.$ The compatibility of the ETNC with respect to
Artin-Verdier/Poitou-Tate duality and the functional equation of
the $L$-function is studied in section 5. A refinement leads to
the formulation of the local $\epsilon$-conjecture  in subsection
5.1.  In order to involve $p$-adic $L$-functions one has to
introduce Selmer groups or better complexes. The necessary
modifications of the $L$-function and the Galois cohomology - in a
way that respects the functional equation -  are described in
section 6. From this the MC in the form of \cite{cfksv} is derived
in subsection 6.2 after a short interlude concerning the new
"localized $K_1$". In the Appendix we collect basic facts about
Galois cohomology on the level of complexes.

{\em Acknowledgements:}  I am very grateful to Takako Fukaya and
Kazuya Kato for providing me with actual  copies of their work and
for answering many  questions. I also would like to thank John
Coates for his interest and Ramdorai Sujatha for reading parts of
the manuscript and  her valuable comments. Finally I am indebted
to David Burns, Pedro Luis del Angel, Matthias Flach, Annette
Huber, Adrian Iovita, Bruno Kahn and Guido Kings for helpful
discussions.

\section{Noncommutative determinants}\label{det}

 The (absolute) TNC measures compares integral structures of Galois
 cohomology with values of complex $L$-functions. For this purpose
 the determinant is the adequate tool as is illustrated by the
 following basic

 \begin{ex}\label{exam}
 Let $T$  be a $\zp$-lattices in a finite dimensional
 $\qp$-vector space $V$ and $f:T\to \mathbb{}T$ a $\zp$-linear map which
 induces an automorphism of $V$. Then the cokernel of $f$ is
 finite with cardinality $ |\det(f)|_p^{-1}$ where $|-|_p$ denotes the
 p-adic valuation normalized as usual: $|p|_p=1/p.$
 \end{ex}

Since the equivariant TNC involves the action of a possibly
non-commutative ring $R$ one needs a determinant formalism over an
arbitrary (associative) ring $R$ (with unit). This can be achieved
by either using virtual objects a la Deligne as Burns and Flach
\cite[\S 2]{bf} do or by  Fukaya and Kato's adhoc construction
\cite[1.2]{fukaya-kato}, both approaches  lead to an equivalent
description.

Let $\mathrm{P}(R)$ denote the category of finitely generated  projective
$R$-modules and $(\mathrm{P}(R),is)$ its subcategory of isomorphisms, i.e.\
with the same objects, but whose morphisms are precisely the
isomorphisms. Then there exists a category $\C_R$ and a functor
\[\d_R:(\mathrm{P}(R),is)\to \C_R\]
which satisfies the following properties:

\begin{enumerate}
  \item[a)] $\C_R$ has an associative and commutative product
  structure $(M,N)\mapsto M\cdot N$ or written just $MN$ with unit object
  $\u_R=\d_R(0)$ and inverses. All objects are of the form
  $\d_R(P)\d_R(Q)^{-1}$ for some $P,Q\in \mathrm{P}(R).$
  \item[b)] all morphisms of $\C_R$ are isomorphisms, $\d_R(P)$ and
  $\d_R(Q)$ are isomorphic if and only if their classes in
  $K_0(R)$ coincide. There is an identification of groups
  $\mathrm{Aut}(\u_R)=K_1(R)$ and $\mathrm{Mor}(M,N)$ is either
  empty or an $K_1(R)$-torsor where $\alpha:\u_R\to \u_R\in K_1(R)$ acts on $\phi:M\to
  N$ as $\alpha\phi:M=\u_R\cdot M\stackrel{\alpha\cdot\phi}{\to}\u_R\cdot N=N.$
  \item[c)] $\d_R$ preserves the "product" structures:  $\d_R(P\oplus Q)=\d_R(P)\cdot\d_R(Q).$
\end{enumerate}

This functor can be naturally extended to complexes. Let
$\mathrm{C}^p(R)$ be the category of bounded complexes in
$\mathrm{P}(R)$ and $(\mathrm{C}^p(R),quasi)$ its subcategory of
quasi-isomorphisms. For $C\in \mathrm{C}^p(R)$  we set
$C^+=\bigoplus_{i\; even}C^i$ and $C^-=\bigoplus_{i\; odd}C^i$ and
define $\d_R(C):=\d_R(C^+)\d_R(C^-)^{-1}$ and thus we obtain a
functor
\[\d_R:(\mathrm{C}^p(R),quasi)\to \C_R\] with the following properties ($C,C',C''\in
\mathrm{C}^p(R)$)

\begin{enumerate}
\item[d)] If $0\to C'\to C\to C''\to 0$ is a short exact sequence
of complexes, then there is a canonical isomorphism
\[\d_R(C)\cong\d_R(C')\d_R(C'')\] which we take as an
identification,
\item[e)] If $C$ is acyclic, then the
quasi-isomorphism $0\to C$ induces a canonical isomorphism
\[\u_R\to\d_R(C).\]
\item[f)] $\d_R(C[r])=\d_R(C)^{(-1)^r}$ where $C[r]$ denotes the $r^\mathrm{th}$ translate of  $C.$
 \item[g)] the functor $\d_R$ factorizes over the image of $\mathrm{C}^p(R)$
in $\mathrm{D}^p(R),$ the category of perfect complexes (as full
triangulated subcategory of the derived category $\mathrm{D}^b(R)$
of the homotopy category of bounded complexes of $R$-modules), and
extends  to  $(\mathrm{D}^p(R),is)$ (uniquely up to unique
isomorphisms) $\footnote{But property d) does not in general
extend  to arbitrary distinguished triangles, thus from a
technical point of view all constructions involving complexes will
have to be made carefully avoiding this problem. We will neglect
this problem but see \cite{bf} for details.}.$
\item[h)] If $C\in \mathrm{D}^p(R)$ has the property that all cohomology
groups $\H^i(C)$ belong again to $\mathrm{D}^p(R),$ then there is
a canonical isomorphism \[\d_R(C)=\prod_i
\d_R(\H^i(C))^{(-1)^i}.\]
\end{enumerate}

Moreover, if $R'$ is another  ring, $Y$ a finitely generated
projective $R'$-module endowed with a structure as right
$R$-module such that the actions of $R$ and $R'$ on $Y$ commute,
then the functor $Y\otimes_R-:\mathrm{P}(R)\to P(R')$ extends to a
commutative diagram

\[\xymatrix{
  (\mathrm{D}^p(R), is) \ar[d]_{Y\otimes_R^\mathbb{L}-} \ar[rr]^{\d_R} & & {\C_R} \ar[d]^{Y\otimes_R-} \\
  (\mathrm{D}^p(R'), is) \ar[rr]^{\d_{R'}} & & {\C_{R'}}   }.\]
In particular, if $R\to R'$ is a ring homomorphism and $C\in \mathrm{D}^p(R),$ we just write $\d_R(C)_{R'}$ for $R'\otimes_R\d_R(C).$

Now let $R^\circ$ be the opposite ring of $R.$ Then the functor
$\mathrm{Hom}_R(-,R)$ induces an anti-equivalence between $\C_R$
and $\C_{R^\circ}$ with quasi-inverse induced  by
$\mathrm{Hom}_{R^\circ}(-,R^\circ);$ both functors will be denoted
by $-^*.$ This extends to a commutative diagram

 \[\xymatrix{
   (\mathrm{D}^p(R), is) \ar[d]_{\mathrm{RHom}_R(-,R)} \ar[rr]^{\d_R}& & {\C_R }\ar[d]^{-^*} \\
   (\mathrm{D}^p(R^\circ),is) \ar[rr]^{\d_{R^\circ}} & & {\C_{R^\circ} } } \]

and similarly for $\mathrm{RHom}_{R^\circ}(-,{R^\circ}).$

For the handling of the determinant functor in practice the
following considerations are quite important:
\begin{rem}\label{inverse}
(i) We have to distinguish at least two inverses of a map
$\phi:\d_R(C)\to \d_R(D)$ with  $C,D\in \mathrm{C}^p(R).$ The
inverse with respect to composition will be denoted by
$\overline{\phi}:\d_R(D)\to\d_R(C).$ But due to the product
structure in $\C_R$ and the identification
$\d_R(C)\cdot\d_R(C)^{-1}=\u_R$ the knowledge of $\phi$ is
equivalent to that of
\[\xymatrix{
  { \u_R}\ar@^{=}[r]&{\d_R(C)\cdot\d_R(C)^{-1}}\ar[rr]^{\phi\cdot \id_{\d_R(C)^{-1}}} & &     {\d_R(D)\cdot \d_R(C)^{-1}}   }
\] or even \[\phi^{-1}:\d_R(C)^{-1}\to \d_R(D)^{-1}\] which is by
definition $\overline{\id_{\d_R(D)^{-1}}\cdot\phi\cdot
\id_{\d_R(C)^{-1}}}$ or in other words
$\phi\cdot\phi^{-1}=\id_{\u_R}. $ In particular, $\phi:\d_R(C)\to
\d_R(C)$ corresponds uniquely to $\phi\cdot
\id_{\d_R(C)^{-1}}:\u_R\to \u_R.$ Thus it can and will be
considered as an element in $K_1(R).$ Note that under this
identification the elements in $K_1(R)$ assigned to each of
$\phi^{-1}$ and $\overline{\phi}$ is the inverse to the element
assigned to $\phi.$ Furthermore, the following relation between
$\circ$ and $\cdot$ is easily verified: Let $\xymatrix{
  A \ar[r]^{\phi} & B   }$ and $\xymatrix{
    B \ar[r]^{\psi} & C   }$ be morphisms in $\C_R.$ Then the composite $\psi\circ\phi$ is the same as  the product $\psi\cdot\phi\cdot\id_{B^{-1}}.$\\
{\bf Convention:} If $\phi:\u\to A$ is a morphism and $B$ an object in $\C_R,$ then we write $\xymatrix{
  B \ar[r]^{\cdot\;\phi} & B\cdot A   }$ for the morphism $\id_B\cdot\phi.$ In particular, any morphism $\xymatrix{
  B \ar[r]^{\phi} & A   }$ can be written as $\xymatrix{
  B \ar[rr]^{\cdot\;(\id_{B^{-1}}\cdot\;\phi)} &&   A   }.$\\
(ii) The determinant of the complex $C=[P_0\stackrel{\phi }{\to}
P_1]$ (in degree $0$ and $1$) with $P_0=P_1=P$ is by definition
$\xymatrix@C=0.5cm{
  { \d_R(C)}\ar@{=}[r]^<(0.3){def} & {\u_R}   }$ and is defined even if $\phi$ is not an
isomorphism (in contrast to $\d_R(\phi)$). But if $\phi$ happens
to be an isomorphism, i.e.\ if $C$ is acyclic, then by e) there is
also a canonical map $\xymatrix@C=0.5cm{
  {  \u_R}\ar[r]^<(0.3){acyc} & {\d_R(C)}   },$ which is in fact  nothing else then
\[\xymatrix@C=0.5cm{ {\u_R}\ar@{=}[r] & {\d_R(P_1)\d_R(P_1)^{-1}
}\ar[rrr]^{\d(\phi)^{-1}\cdot \id_{\d(P_1)^{-1}}} &&&
{\d_R(P_0)\d_R(P_1)^{-1}}   \ar@{=}[r] & {\d_R(C)} }\] (and  which
depends in contrast to the first identification on $\phi$). Hence,
the composite $\xymatrix@C=0.5cm{
  {  \u_R}\ar[r]^<(0.3){acyc} & {\d_R(C)} \ar@{=}[r]^<(0.4){def} & {\u_R}
  }$ corresponds to $\d_R(\phi)^{-1}\in K_1(R)$ according to the first
remark.
 In order to distinguish the above identifications between $\u_R$ and $\d_R(C)$     we also say that $C$ is {\em trivialized by
the identity } when we refer to $\xymatrix@C=0.5cm{
  { \d_R(C)}\ar@{=}[r]^<(0.4){def} & {\u_R}   }$ (or its inverse with
  respect to composition). For $\phi=\id_P$ both identifications
   agree obviously.
 \end{rem}

We end this section by considering the example where $R=K$ is a
field and $V$ a finite dimensional vector space over $K.$ Then,
according to \cite[1.2.4]{fukaya-kato}, $\d_K(V)$ can be
identified with the highest exterior product $\bigwedge^{top}V$ of
$V$ and for an automorphism $\phi: V\to V$ the determinant
$\d_K(\phi)\in K^\times=K_1(K)$ can be identified with the usual
determinant $\det_K(\phi).$ In particular, we identify $\d_K=K$
with canonical basis $1.$ Then a map $\xymatrix@C=0.5cm{ {\u_K}
\ar[r]^{\psi} & {\u_K}  }$ corresponds uniquely to the value
$\psi(1)\in K^\times.$

\begin{rem}\label{finitemodules}
Note that every {\em finite} $\zp$-module $A$ possesses a free
resolution $C$   as in Remark \ref{inverse} (ii), i.e.\
$\d_\zp(A)\cong\d_\zp(C)^{-1}=\u_\zp.$ Taking into account the
above and Example \ref{exam}  we see that modulo $\z_p^\times$ the
composite $\xymatrix@C=0.5cm{
  {  \u_\qp}\ar[r]^<(0.2){acyc} & {\d_\zp(C)_\qp} \ar@{=}[r]^<(0.4){def} & {\u_\qp}
  }$ corresponds to the cardinality $|A|\in\Q_p^\times.$
\end{rem}

\section{$K$-Motives over $\Q$}

In this survey we will be mainly interested in the Tamagawa Number
Conjecture and Iwasawa theory for the motive $M=h^1(E)(1)$ of an
elliptic curve $E$ or the slightly more general $M=h^1(A)(1)$ of
an abelian variety $A$  defined over $\Q.$ But as it will be
important to consider certain twists of $M$ we also recall basic
facts on the Tate motive $\Q(1)$ and Artin motives. We shall
simply view motives in the naive sense, as being defined by a
collection of realizations satisfying certain axioms, together
with their motivic cohomology groups. The archetypical motive is
$h^i(X)$ for a smooth projective variety $X$ over $\Q$ with its
obvious {\'e}tale cohomology
$\H^i_{\acute{e}t}(X\times_\Q\overline{\Q},\Q_l)$, singular
cohomology $\H^i(X(\mathbb{C}),\Q)$ and de Rham cohomology
$\H^i_{dR}(X/\Q),$ their additional structures and comparison
isomorphisms. More general, let $K$ be a finite extension of $\Q.$
A $K$-motive $M$ over $\Q,$ i.e.\ a motive over $\Q$ with an
action of $K,$ will be given by the following data, which for
$M=h^n(X)_K$ arise by tensoring the above cohomology groups by $K$
over $\Q :$

\subsection{The $l$-adic realization  $M_l$ of $M$ (for every prime number $l$)}

For   a place $\lambda$ of $K$  lying above $l$ and denote by
$K_\lambda$ the completion of $K$ with respect to $\lambda.$ Then
$M_\lambda$ is a continuous finite dimensional $K_\lambda$-linear
representation of the absolute Galois group $G_\Q$ of $\Q.$ We put
$K_l:=K\otimes_\Q \Q_l=\prod_{\lambda |l} K_\lambda$ and we denote
by $M_l$ the free $K_l$-module $\prod_{\lambda |l} M_\lambda.$

\subsection{The Betti realization $M_B$ of $M$}

Attached  to $M$ is a finite dimensional $K$-vector space $M_B$
which carries an action of complex conjugation $\iota$ and a
$\Q$-Hodge structure  $M_B\otimes_\Q \bbC \cong \bigoplus
\h^{i,j}$ (over $\R$) with $\iota\h^{i,j}=\h^{j,i}$ where
$\h^{i,j}$ are free $K_\bbC:=K\otimes_\Q \bbC\cong
\bbC^{\Sigma_K}$-modules and where $\Sigma_K$ denotes the set of
all embeddings $K\to \bbC.$ E.g.\   the motive $M=h^n(X)$ is pure
of weight $w(M)=n,$ i.e.\ $\h^{i,j}=0$ if $i+j\neq n.$

\subsection{The de Rham realization $M_{dR}$ of $M$}

$M_{dR}$ is  a finite dimensional $K$-vector space with a
decreasing exhaustive filtration $M^k_{dR},$ $k \in \z$. The
quotient $t_M=M_{dR}/M^0_{dR}$ is called the {\em tangent space}
of $M.$

\subsection{Comparison between $M_B$ and $M_l$}

For each prime number $l$ there is an isomorphism of $K_l$-modules
\begin{equation}
 \xymatrix{
  K_l\otimes_K M_B \ar[r]^<(0.4){g_l}_<(0.4){\cong} & M_l   }
\end{equation}
which respects the action of complex conjugation, in  particular
it induces  canonical isomorphisms
\begin{equation}\label{B-l}
g_\lambda^+:K_\lambda\otimes_K M_B^+ \cong M_\lambda^+     \mbox{
and } g_l^+:K_l\otimes_K M_B^+ \cong M_l^+.
\end{equation}
Here and in what follows, for any commutative ring $R$ and
$R[G(\bbC/\R)]$-module $X$ we denote by $X^+$ and $X^-$ the
$R$-submodule of $X$ on which $\iota$ acts by $+1$ and $-1,$
respectively.

\subsection{Comparison between $M_B$ and $M_{dR}$} There is a
$G(\bbC/\R)$-invariant isomorphism of $K_\bbC$-modules
\begin{equation}\label{B-dR}
 \xymatrix{
{\bbC\otimes_\Q M_B }\ar[r]^{g_\infty}_{\cong} & {\bbC\otimes_\Q
M_{dR}} }
\end{equation}
(on the left hand side $\iota$ acts diagonally while on the right
hand side only on $\bbC$) such that for all $k\in \z$

\[g_\infty(\bigoplus_{i\geq k }\h^{i,j}(M))\cong  \bbC\otimes_\Q M_{dR}^k.\]

This induces an isomorphism
\begin{equation}
\label{B+-dR} (\bbC\otimes_\Q M_B)^+\cong \R\otimes_\Q M_{dR}
\end{equation}
and the period map
\begin{equation}\label{alphaM}
 \xymatrix{
{\R\otimes_\Q M_B^+} \ar[r]^{\alpha_M} & {\R\otimes_\Q
t_M} }
\end{equation}

We say that $M$ is {\em critical} if this happens to be an
isomorphism\footnote{By \cite[lem.\ 3]{coates91} $M$ is critical
if and only if one of the following equivalent conditions holds:
a) both infinite Euler factors $L_\infty(M,s)$ and
$L_\infty(M^*(1),-s)$ (see section \ref{functional}) are
holomorphic at $s=0,$ b) if $j<k$ and $\h^{j,k}\neq \{0\}$ then
$j<0$ and $k\geq 0,$ and, in addition, if $\h^{k,k}\neq \{0\},$
then $\iota$ acts on this space as $+1$ if $k<0$ and by $-1$ if
$k\geq 0.$ See also \cite[lem.\ 2.3]{coates-perrin-riou} for
another criterion.}.

\subsection{Comparison between $M_p$ and $M_{dR}$}

Let $B_{dR}$ be the filtered field of de Rham periods with respect
to $\overline{\qp}/\qp,$ which is endowed with a continuous action
of the absolute Galois group $G_\qp$ of $\qp,$ and set as usual
$D_{dR}(V)=(B_{dR}\otimes_\qp V)^{G_\qp}$ for a finite-dimensional
$\qp$-vector space $V$ endowed with a continuous action  of
$G_\qp.$ The  (decreasing) filtration $B_{dR}^i$ of  $B_{dr}$
induces a filtration  $D_{dR}^i=(B_{dR}^i\otimes_\qp V)^{G_\qp}$
of $D_{dR}.$ Then there is an $G_{\qp}$-invariant isomorphism of
filtered $K_p\otimes_\qp B_{dR}$-modules
\begin{equation}
 \xymatrix{
 {B_{dR}\otimes_\qp  M_p} \ar[r]^<(0.2){g_{dR}}_<(0.2){\cong} &{B_{dR}\otimes_\Q M_{dR}}}
\end{equation}
which induces an  isomorphism of filtered $K_p$-modules by taking
$G_\Q$-invariants
\begin{equation}
 \xymatrix{
{  D_{dR}(M_p)} \ar[r]^{g_{dR}}_{\cong} & { K_p\otimes_K M_{dR}},
}
\end{equation}
an isomorphism of $K_p$-modules
\begin{equation}\label{dR-p}
 \xymatrix{
{ t(M_p):= D_{dR}(M_p)/D_{dR}^0(M_p)}
\ar[r]^<(0.2){g_{dR}^t}_<(0.2){\cong} & { K_p\otimes_K t_M} }
\end{equation}
and, for each place $\lambda$ of $K$ over $p,$ an isomorphism of
$K_\lambda$-vector spaces
\begin{equation}
 \xymatrix{
{  t(M_\lambda):=D_{dR}(M_\lambda)/D_{dR}^0(M_\lambda)}
\ar[r]^<(0.2){g_{dR}^t}_<(0.2){\cong} & { K_\lambda\otimes_K t_M.}
}
\end{equation}

The tensor product $M\otimes_K N$ of two $K$-motives  is given by
the data which arises from the tensor products of all realizations
and their additional structures. Similar the dual $M^*$ of the
$K$-motive $M$ is given by the duals of the corresponding
realizations. In particular, we denote by $M(n),$ $n\in \z,$ the
twist of $M$ by the $|n|$-fold tensor product
$\Q(n)=\Q(1)^{\otimes n}$ of the Tate motive if $n\geq 0$ and of
its dual $\Q(-1)=\Q(1)^*$ if $n<0.$ For the motive $M=h^i(X)(j)$
where the dimension of $X$ is $d,$ Poincar{\'e} duality gives a
perfect pairing
\[h^i(X)(j)\times h^{2d-i}(X)(d-j)\to h^{2d}(X)(d)\cong \Q\]
which identifies  $M^*$ with $h^{2d-i}(X)(d-j).$ Here
$\Q=h^0(spec(\q))(0)$ denotes the trivial $\Q$-motive.

\begin{ex}
A) The Tate motive $\Q(1)=h^2(\mathbb{P}^1)^*$ should be thought
of as $h_1(\mathbb{G}_m)$ even though the multiplicative group
$\mathbb{G}_m$ is not proper. Its $l$-adic realisation is the
usual Tate module $\ql(1)$ on which $G_\Q$ acts via the cyclotomic
character $\chi_l: G_\Q\to \mathbb{Z}_l^\times.$ The action of
complex conjugation on $\Q(1)_B=\Q$ is by $-1,$ its Hodge
structure is pure of weight $w(M)=-2$ and given by $\h^{-1,-1}.$
The filtration $\Q(1)_{dR}^k$ of $\Q(1)_{dR}=\Q$ is either $\Q$ or
$0,$ according as $k\leq -1$ or $k>-1,$ in particular we have
$t_{\Q(1)}=\Q.$ Finally, $g_\infty$ sends $1\otimes 1$ to $2\pi
i\otimes 1$ while $g_{dR}$ sends $1\otimes 1$ to $t\otimes 1,$
where $t="2\pi i"$ is the $p$-adic period analogous to $2\pi i.$

B) For the $\Q$-motive $M=h^1(A)(1)$ of an abelian variety $A$
over $\Q$ we have
$M_l=H_{\acute{e}t}^1(A_{\bar{\Q}},\ql(1))=\Hom_\ql(V_lA,\ql(1))\cong
V_l(A^\vee)$ via the Weil pairing. More generally, the Poincare
bundle on $A\times A^\vee$ induces  isomorphisms $M^*(1)\cong
h^1(A^\vee)(1)$ and $M\cong h^1(A^\vee)^*,$ while by fixing a
(very) ample symmetric line bundle on $A,$  whose existence is
granted by \cite[cor.\ 7.2]{mil-abvar}, it is sometimes convenient
to identify $M$ with $h_1(A):=h^1(A)^* $ using the hard Lefschetz
theorem (\cite[1.15,thm.\ 5.2 (iii)]{scholl-motives}, see also
\cite{kleiman}) (but in general better to work with the dual
abelian variety $A^\vee).$ Then $M_l$ can be identified with
$V_l(A),$ while $M_B=H^1(A(\bbC),\Q)(1)$ can be identified with
$H_1(A(\bbC),\Q),$ the Hodge-decomposition (pure of weight $-1$)
is given by $\h^{0,-1}=H^0(A(\bbC),\Omega_A^1)(\cong
\Hom_\bbC(\H^1(A(\bbC),\Omega_A^0),\bbC))$ and
$\h^{-1,0}=\H^1(A(\bbC),\Omega^0_A)(\cong\Hom_\bbC(\Omega^1(A),\bbC)).$
Furthermore, we have $M_{dR}^{-1}=M_{dR},$ $M^0_{dR}$ is the image
 of $\Omega^1_{A/\Q}(A)(\cong
\H^1(A,\Omega^0_{A/\Q})^*)$ and $M^1_{dR}=0.$ In particular, $t_M=
\H^1(A,\Omega^0_{A/\Q})=Lie(A^\vee)$ (e.g.\ \cite[thm.
5.11]{klei}) the Lie-algebra of $A^\vee,$ can be identified with
$t_{h_1(A)}=\Hom_\Q(\Omega^1_{A/\Q}(A),\Q)=Lie(A).$ The map
$\alpha_M$  for the motive $M=h_1(A),$ which is in fact an
isomorphism, is induced by sending a $1$-cycle $\gamma\in
\H_1(A(\bbC),\Q)^+ $ to $\int_\gamma\in
\Hom_\Q(\Omega^1_{A/\Q}(A),\bbR)=Lie(A)_\bbR$ which sends a
$1$-form $\omega$ to $\int_\gamma\omega\in \bbR.$

C) Artin motives $[\rho]$ (with coefficients in  a finite
extension $K$ of $\Q$) are direct summands of the $K$-motive
 $h^0(spec(F))\otimes_\Q K$ but can also be identified with the
 category of  finite-dimensional $K$-vector spaces $V$ with an
 action by $G_\Q,$ i.e.\   representations $\rho:G_\Q\to Aut_K(V)$
 with finite image. We write $[\rho]$ for the corresponding motive
 and have $[\rho]_l=V\otimes_K K_l$ with $G_\Q$ acting just on
 $V,$ $[\rho]_B=V$ with Hodge-Structure pure of Type $(0,0)$ and
 $[\rho]_{dR}=(V\otimes_\Q \bar{\Q})^{G_\Q},$ where $G_\Q$ acts
 diagonally. Since $[\rho]_{dR}^k$ is either $[\rho]_{dR}$ or $0$ according as $k\leq 0$ of $k>0,$ we have $t_{[\rho]}=0.$ The inverse of $g_\infty$ is induced by the natural
 inclusion $(V\otimes_\Q \bar{\Q})^{G_\Q}\subseteq  V\otimes_\Q
 \bar{\Q}.$ E.g.\ if $\psi$ denotes a Dirichlet character of
 conductor $f$ considered via $(\z/f\z)^*\cong G(\Q(\zeta_f)/\Q)$
 as character $G\to K^\times$ where $K=\Q(\zeta_{\varphi(f))})$ and $\varphi$ denotes the Euler $\varphi$-function, then we
 obtain a basis of $[\psi]_{dR}$ over $K$ by the Gauss sum
 \[\sum_{1\leq n <f,(n,f)=1} \psi(n)\otimes e^{-2\pi i n/f}\;\;\in (K(\psi)\otimes_\Q \bar{\Q})^{G_\Q},\]
 where $K(\psi)$ denotes the $1$-dimensional $K$-vector space on which $G_\Q$ acts via $\psi.$

 Of course,  $h^0(spec(F))\otimes_\Q K$ corresponds to the regular
 representation of $G(F/\Q)$ on $K[G(F/\Q)]$ considered as
 representation of $G_\Q.$

 Other examples arise by taking symmetric products or tensor
 products of the above examples. In particular, we will be
 concerned with the  motives

D) $[\rho]\otimes h^1(A)(1),$ where $\rho$ runs through all Artin
representations.\\
 E) Finally, the motive $M(f)$ of a modular form is a prominent
 example, see \cite[\S 7]{deligne-L} and \cite{scholl}.
\end{ex}

\subsection{Motivic cohomology}
The motivic cohomology $K$-vector spaces \linebreak $\H_f^0(M):=\H^0(M)$ and
$\H^1_f(M)$ may be defined by algebraic $K$-theory or motivic
cohomology a la Voevodsky. They are conjectured to be finite
dimensional. Instead of a general definition we just describe them
in our standard examples.

\begin{ex}\label{ex-mot-coh}
A) For the Tate motive we have $\H_f^0(\Q(1))=\H_f^1(\Q(1))=0$ and
for its Kummer dual $\H_f^0(\Q)=\Q$ while $\H_f^1(\Q)=0.$

B) If $M=h^1(A)(1)$  for an abelian variety $A$ over $\Q$ one has $\H^0_f(M)=0$ and $\H^1_f(M)=A^\vee(\Q)\otimes_\z\Q.$

C) For $M=h^0(spec(F))$  we have $\H_f^0(M)=\Q$ and $\H_f^1(M)=0$
while for $M^*(1)=h^0(spec(F))(1)$ one has $\H_f^0(M^*(1))=0$ and
$\H_f^1(M^*(1))= \O_F^\times\otimes_\z\Q.$ More general, for an
$K$-Artin motive $[\rho]$ one has $\H^0_f([\rho])=K^n,$ where $n$
is the multiplicity with which $\Q$ occurs in $[\rho].$

Unfortunately the functor $\H^i_f$ does not behave well with
tensor  products, i.e.\ in general one cannot derive
$\H^*_f([\rho]\otimes h^1(A)(1))$ from $\H^*_1([\rho])$ and
$\H^*_1(h^1(A)(1))$ (e.g.\ in form of a K\"unneth formula).
\end{ex}

\section{The Tamagawa Number Conjecture - absolute version}

In \cite{BK} Bloch and Kato formulated a vast generalization of
the analytic class number formula and the BSD-conjecture. While
the conjecture of Deligne and Beilinson links the order of
vanishing of the $L$-function attached to a motive $M$ to its
motivic cohomology and claims rationality of special $L$-values or
more general  leading coefficients (up to periods and regulators)
the Tamagawa number conjecture by Bloch and Kato predicts the
precise $L$-value in terms of Galois cohomology (assuming the
conjecture of Deligne-Beilinson).

Later, Fontaine and Perrin-Riou \cite{fp} found an equivalent
formulation using (commutative) determinants instead of (Tamagawa)
measures\footnote{The name comes from an analogy with the theory
of algebraic groups, see \cite{BK}.}. In this section we follow
closely their approach.

Let us first recall the definition of the complex $L$-function
attached to a $K$-motive $M.$ We fix a place $\lambda$ of $K$
lying over $l$ and an embedding $K\to \bbC.$ For every prime $p$
take a prime $l\neq p$ and set
\[P_p(M_\lambda,X)=\mathrm{det}_{K_\lambda}(1-\varphi_p X|(M_\lambda)^{I_p})\in K_\lambda[X],\]

where $\varphi_p$ denotes the geometric Frobenius automorphism of
$p$ in $G_{\qp}/I_p$ and $I_p$ is the inertia subgroup of $p$ in
$G_\qp\subseteq G_\Q.$ It is conjectured that $P_p(X)$ belongs to
$K[X]$ and is independent  of the choices of $l$ and $\lambda.$
For example this is known by the work of Deligne proving the Weil
conjectures for $M=h^i(X)$ for places $p$ where $X$ has good
reduction; by the compatibility of the system of $l$-adic
realisations for abelian varieties \cite[rem.
2.4.6(ii)]{fontaine-ss} and Artin motives it is also clear for our
examples A)-D). Then we have the $L$-function of $M$ as Euler
product
\[L_K(M,s)=\prod_p P_p(M_\lambda,p^{-s})^{-1},\]
defined and analytic for $\Re(s)$ large enough.

\begin{ex}
A) The $L$-function $L_\Q(\Q(1),s-1)$ of the Tate motive is just the Riemann zeta function $\zeta(s).$ In general, one has $L_K(M(n),s)=L_K(M,s+n)$ for any $K$-motive $M$ and any integer $n.$\\
B) If $M=h^1(A)(1)$  for an abelian variety $A$ over $\Q,$ then $L(M,s-1)$ is the classical Hasse-Weil $L$-function of $A^\vee,$ which coincides with that for $A$ because $A$ and $A^\vee$ are isogenous.\\
C) $L_K([\rho],s)$ coincides with the usual Artin $L$-function of $\rho,$ in particular we retrieve the Dedekind zeta-function $\zeta_F(s)$ as $L_\Q(h^0(spec(F),s).$  \\
D) The $L$-functions $L_K([\rho]\otimes h^1(A)(1),s)$ will play a
crucial role for the interpolation property of the $p$-adic
$L$-function.
\end{ex}

Also, the meromorphic continuation to the whole plane $\bbC$ is
part of the conjectural framework. The Taylor expansion
\[L_K(M,s)=L^*_K(M) s^{r(M)} + \ldots\]

defines the leading coefficient $L^*_K(M) \in \bbC^\times,$ which
can be shown to belong to $\R^\times$ actually, and the order of
vanishing $r(M)\in \z$ of $L_K(M,s)$ at $s=0.$ The aim of the
conjectures to be formulated now is to express $L^*(M)$ and $r(M)$
in terms of motivic and Galois cohomology.

\begin{conj}[Order of Vanishing; Deligne-Beilinson]\label{order}
 \[r(M)=\dim_K \H^1_f(M^*(1)) -\dim_K \H_f^0(M^*(1))\]
\end{conj}

According to the remark in \cite{flach-survey} the duals of
$\H^i_f(M^*(1))$ should be considered as  "motivic cohomology with
compact support $\H_c^{2-i}(M)$" and thus $r(M)$ is just their
Euler characteristic. This explains why the Kummer duals $M^*(1)$
are involved here.

The link between the complex world, where the values $L^*(M)$
live, and the $p$-adic world, where the Galois cohomology lives,
is formed by the {\em fundamental line} in $\C_K$ following the
formulation of Fontaine and Perrin-Riou \cite{fp}:

\begin{eqnarray*}
\Delta_K(M):&=&\d_K(\H^0_f(M))^{-1}\d_K(\H^1_f(M))\d_K(\H^0_f(M^*(1))^*)
\d_K(\H^1_f(M^*(1))^*)^{-1}\\ && \d_K(M_B^+)\d_K(t_M)^{-1}.
\end{eqnarray*}

The relation of $\Delta_K(M)$ with the Betti and de Rham
realization of $M$  is given by the following

\begin{conj}[Fontaine/Perrin-Riou]\label{motinfty}
There exist an exact sequence of  $K_{\mathbb{R}}:=\R\otimes_\Q
K$-modules
\[\begin{split}
&\xymatrix{
  0 \ar[r] & { \H^0_f(M)_\mathbb{R}\ar[r]^{c}}
  & {\ker(\alpha_M)}  \ar[r]^<(0.3){r_B^*}
  & {(\H^1_f(M^*(1))_\mathbb{R})^* }  & }\\
  &\xymatrix{ {\phantom{0}}   {\ar[r]^<(0.3){h}} & {\H^1_f(M)_\mathbb{R}} \ar[r]^{r_B} & {\coker(\alpha_M)} \ar[r]^{c^*} & {(\H^0_f(M^*(1))_{\mathbb{R}})^*\ar[r]  }& {0   }}
\end{split}  \]
where by $-_\mathbb{R}$ we denote the base change  from $\Q$ to
$\R$ (respectively from $K$ to $K_\mathbb{R}$), $c$ is the cycle
class map into singular cohomology, $r_B$ is the Beilinson
regulator map and (if both $\H^1_f(M)$ and $\H^1_f(M^*(1))$ are
nonzero so that $M$ is of weight $-1, $ then) $h$ is a height
pairing.
\end{conj}

\begin{ex} A),C) For the motive $M=h^0(spec(F))$  the above exact
sequence is just the $\R$-dual of the following
\[\xymatrix@C=0.5cm{
  0 \ar[r] & {\O_F^\times\otimes_\z\R} \ar[rr]^{r} && {\R^{r_1}\times\bbC^{r_2}} \ar[rr]^{ \Sigma} && {\R} \ar[r] & 0
  },\] where $r$ is the Dirichlet(=Borel) regulator map
  (see \cite{bur} for a comparison of the Beilinson and Borel
  regulator map)
  and $r_1$ and $r_2$ denote the number of real and complex places
  of $F,$ respectively.\\
B) The Neron height pairing (see \cite{bloch} and the references
there) \[<,>: A^\vee(\Q)\otimes_\z\R \times A(\Q)\otimes_\z\R \to
\R \] induces a homomorphism $A^\vee(\Q)\otimes_\z\R\to
\Hom_\z(A(\Q),\R),$ the inverse of which gives the exact sequence
for the motive $M=h^1(A)(1).$
\end{ex}

We assume this conjecture. Using property e) and change of rings
of the functor $\d,$ it induces a canonical isomorphism
(period-regulator map)
\begin{equation}
\vartheta_\infty: K_\bbR\otimes_K\Delta_K(M)\cong \u_{K_\bbR}.
\end{equation}

\begin{conj}[Rationality; Deligne-Beilinson]\label{rationality}
 There is a unique isomorphism
 \[\zeta_K(M): \u_K \to \Delta_K(M)\]
 such that for every embedding $K \to \bbC$ we have
 \[ \xymatrix{
   L^*_K(M):\u_\bbC\ar[rr]^{\zeta_K(M)_\bbC} &   & {\Delta_K(M)_\bbC} \ar[rr]^{(\vartheta_\infty)_\bbC} &  & {\u_\bbC}   }\]
\end{conj}

In other words, the preimage
$\overline{\vartheta_\infty}(L^*_K(M))$ of $L^*_K(M)$ with respect
to $\vartheta_\infty$ generates the $K$-vector space $\Delta_K(M)$
if the determinant functor is identified with the  highest
exterior product. Thus  up to a period and a regulator (the
determinant of $\vartheta_\infty$ with respect to a $K$-rational
basis) the value $L^*_K(M)$ belongs to $K.$

The rationality enables us to relate $L^*_K(M)$ to the $p$-adic
world which we will describe now.

Let $S$ be a finite set of places of $\Q$ containing $p,\infty$
and the places of bad reduction of $M,$ then
$U:=spec(\z[\frac{1}{S}])$ is an open dense subset of $spec(\z ).$
Then we have complexes $\mathrm{R\Gamma}_c(U,M_p),$ $
\mathrm{R\Gamma}_f(\Q,M_p)$ and $\mathrm{R\Gamma}_f(\Q_v,M_p)$
 calculating the (global) cohomology $\H^i_c(U,M_p)$ with compact
 support, the finite part of global and local cohomology,
 $\H^i_f(\Q,M_p)$ and $\H^i_f(\Q_v,M_p),$ respectively, see
 \ref{appendix}. These complexes fit into a distinguished
 triangle (see  \eqref{trian-c-f})
 \begin{equation}\label{triangle_c_f_f}
  \xymatrix{
  {\mathrm{R\Gamma}_c(U,M_p)} \ar[r]  & {\mathrm{R\Gamma}_f(\Q,M_p)} \ar[r]  & {\bigoplus_{v\in S}\mathrm{R\Gamma}_f(\Q_v,M_p) } \ar[r]  &    . }
 \end{equation}
On the other hand motivic cohomology specializes to the finite
parts of global Galois cohomology:

\begin{conj}\label{motl}
 There are natural isomorphisms  $\H^0_f(M)_{\Q_l}\cong
 \H_f^0(\Q,M_l)$ (cycle class maps) and $\H^1_f(M)_{\Q_l}\cong
 \H_f^1(\Q,M_l)$ (Chern class maps).
\end{conj}

Hence, as there is a duality $\H^i_f(\Q,M_l)\cong
\H_f^{3-i}(\Q,M_l^*(1))^*$ for all $i,$ this conjecture determines
all cohomology groups $\H_f^i(\Q,M_l).$


Using   properties d), g) and change of rings of the determinant
functor, the conjecture \ref{motl} for $l=p$, the canonical
isomorphisms (see appendix \ref{eta})
\begin{eqnarray}
 \eta_p(M_p):& \u_{K_p}& \to
 \d_{K_p}(\mathrm{R\Gamma}_f(\Q_p,M_p))\cdot
 \d_{K_p}(t(M_p)),\\
\eta_l(M_p):& \u_{K_p} &\to
 \d_{K_p}(\mathrm{R\Gamma}_f(\Q_l,M_p)),
\end{eqnarray}
the comparison isomorphisms \eqref{B-l} and \eqref{dR-p} as well
as \eqref{triangle_c_f_f}, we obtain a canonical isomorphism
($p$-adic period-regulator map)
\begin{equation}\label{vartheta}
\vartheta_p:  \Delta_K(M)_{K_p}\cong
\d_{K_p}({\mathrm{R\Gamma}_c(U,M_p)})^{-1} ,
\end{equation}
 which induces for any place $\lambda$ above $p$
\begin{equation}
\vartheta_\lambda:  \Delta_K(M)_{K_\lambda}\cong
\d_{K_\lambda}({\mathrm{R\Gamma}_c(U,M_\lambda)})^{-1} .
\end{equation}

Now let $T_\lambda$ be a Galois stable $\O_\lambda$-lattice of
$M_\lambda$ and $\mathrm{R\Gamma}_c(U,T_\lambda)$ its Galois
cohomology with compact support, see section \ref{appendix}. Here,
$\O_\lambda$ denotes the  valuation ring of $K_\lambda.$ Note that
by Artin-Verdier/Poitou-Tate  duality (see \ref{AVPT}) the
"cohomology" $\mathrm{R\Gamma}_c(U,T_\lambda)$ with compact
support can also be replaced by the complex $\r(U,
T_\lambda^*(1))^*\oplus (T_\lambda^*(1))^+$ where $\r(U,
T_\lambda^*(1))$ calculates as usual the global Galois cohomology
with restricted ramification.

The following conjecture, for every prime $p,$  gives a precise
description of the special $L$-value $L^*_K(M)\in \bbR^\times$ up
to $\O_K^\times,$ i.e.\ up to sign if $K=\Q,$ where $O_K$ denotes
the ring of integers in $K$:

\begin{conj}[Integrality; Bloch/Kato, Fontaine/Perrin-Riou]\label{integrality}
 Assume conjecture \ref{rationality}. Then for every place $\lambda$ above $p$ there exist a
 (unique)  isomorphism
 \[\zeta_{\O_\lambda}(T_\lambda): \u_{\O_\lambda} \to \d_{\O_\lambda}(\mathrm{R\Gamma}_c(U,T_\lambda))^{-1}\]
 which induces via $K_\lambda\otimes_{\O_\lambda} -$ the following
 map
 \be\label{zetaintegral}\xymatrix{
   {\zeta_{\O_\lambda}(T_\lambda)_{K_\lambda}:\u_{K_\lambda}} \ar[rr]^{\zeta_K(M)_{K_\lambda}} &
    & {\Delta_K(M)_{K_\lambda}} \ar[rr]^{\vartheta_\lambda} &   & {\d_{K_\lambda}({\mathrm{R\Gamma}_c(U,M_\lambda)})^{-1}.}
       }\ee
\end{conj}

If we identify again the determinant functor with the highest
exterior product, this conjecture can be rephrased as follows:
$\vartheta_\lambda\overline{\vartheta_\infty}(L_K^*(M))$ generates
the $\O_\lambda$-lattice \linebreak
$\d_{K_\lambda}({\mathrm{R\Gamma}_c(U,T_\lambda)})^{-1}$ of
$\d_{K_\lambda}({\mathrm{R\Gamma}_c(U,M_\lambda)})^{-1}.$  In
other words, this generator is determined up to a unit in
$\O_\lambda.$

It can be shown that this conjecture is independent of the choice
of $S$ and $T_\lambda.$

\begin{ex}
(Analytic class number formula) For the motive $M=h^0(spec(F))$ we
have  that $r(M)=r_1+r_2-1$ if $F\otimes_\Q
\bbR\cong\bbR^{r_1}\times \bbC^{r_2}$ and that
$L^*(M)=\frac{-|Cl(\O_F)|R}{|\mu(F)|}$ for the unit regulator $R.$
Thus Conjectures \ref{order},\ref{rationality} and
\ref{integrality} are theorems in this case!
\end{ex}

For other known cases of these conjectures we refer the reader to the excellent survey \cite{flach-survey}, where in particular the results of Burns-Greither \cite{burns-greither2003} and Huber-Kings \cite{hu-ki2} are discussed.

Another example will be discussed in the following section.

\subsection{Equivalence to classical formulation of
BSD}\label{BSD}
I am very grateful to Matthias Flach for some advice concerning
this section, in which  we assume $p\neq 2$  for simplicity. In
order to see that the above conjectures for the motive
$M=h^1(A)(1)$ of an abelian variety $A$ are equivalent to the
classical formulation involving all the arithmetic invariants of
$A$ one has to consider also integral structures for the finite
parts of global and local Galois cohomology. For $T_p$ we take the
Tate-module $T_p(A^\vee)$ of $A^\vee.$ In particular one can
define perfect complexes of $\zp$-modules $\r_f(\Q,T_p)$ and
$\r_f(\Q_v,T_p)$ such that the analogue of \eqref{triangle_c_f_f}
holds, see \cite[{\S}1.5]{bf3}. We just state some results concerning
their cohomology groups $\H_f^i$ in the following

\begin{prop}[{\cite[(1.35)-(1.37)]{bf3}}]\label{resultsabelianvar}
(a)(global) If the Tate-Shafarevich group $\Sha(A/\Q)$ is finite,
then one has
\begin{align}
\H^0_f(\Q,T_p)&=0 & \H^3_f(\Q,T_p)&\cong\Hom_\z(A(\Q)_{tors},\qp/\zp) \\
\H^1_f(\Q,T_p)&\cong A^\vee(\Q)\otimes_\z \zp & \H^i_f(\Q,T_p)&=0
\mbox{ for } i\neq 0,1,2,3
\end{align}
 and an exact sequence of $\zp$-modules
 \be \xymatrix@C=0.5cm{
   0 \ar[r] & {\Sha(A/\Q)(p)} \ar[rr]^{ } && {\H^2_f(\Q,T_p)} \ar[rr]^{ } && {\Hom_\z(A(\Q),\zp)} \ar[r] & 0.}\ee
(b)(local) For all primes $l$ one has
\begin{align}
\H^0_f(\ql,T_p)&=0 & \H^1_f(\ql,T_p)&\cong A^\vee(\ql)^{\wedge p}
& \H^i_f(\ql,T_p)&=0 \mbox{ for } i\neq 0,1,
\end{align}
where $A^\vee(\ql)^{\wedge p}$ denotes the $p$-completion of
$A^\vee(\ql).$
\end{prop}

Note that one has $\H^i_f(-,T_p)\otimes_\zp \qp\cong
\H^i_f(-,M_p)$ for the local and global versions.

Recall from Example \ref{ex-mot-coh} B) that we have \be
\Delta_\q(M)=\d_\q(A^\vee(\q)\otimes_\z\q))\d_\q(\Hom_\z(A(\q),\q))^{-1}
\d_\q(H_1(A^\vee(\bbC),\Q)^+)\d_\q(Lie(A^\vee))^{-1}\ee

In order to define  the period and regulator we have to choose
bases: we first fix $P_1^\vee,\ldots ,P_r^\vee\in A^\vee(\q)$
(respectively $P_1,\ldots ,P_r\in A(\q)$),
$r=\rk_\z(A^\vee(\q))=\rk_\z(A(\q)),$ such that setting
$T_{A^\vee}:=\bigoplus \z P_i^\vee$ (respectively $T_A^d:=
\bigoplus \z P_i^d\subseteq\Hom_\z(A(\q),\z)$, where $P_i^d$
denotes the obvious dual "basis") we obtain
 \begin{align}
A^\vee(\q)&\cong  A^\vee(\q)_{tor}\oplus T_{A^\vee} &
\Hom_\z(A(\q),\z)&\cong T_A^d.
 \end{align}

 Similarly we fix a $\z$-basis $\gamma^+=(\gamma_1^+,\ldots, \gamma_{d_+}^+)$ of
 $T_B^+:=H_1(A^\vee(\bbC),\z)^+$ and $\z$-basis
 $\delta=(\delta_1,\ldots,\delta_{d_+})$ of the $\z$-lattice $Lie_\z(A^\vee):=\Hom_\z(\Omega^1_{\B/\z}(\B),\z)$ of $Lie(A^\vee),$
 respectively. Here $\B/\z $ denotes the (smooth, but not proper) Neron
 model of $A^\vee$ over $\z.$  Thus we obtain an integral structure of $\Delta_\Q(M):$
 \be \Delta_\z(M):=
 \d_\z(T_{A^\vee})\d_\z(T_A^d)^{-1}\d_\z(T_B^+)\d_\z(Lie_\z(A^\vee))^{-1}\ee
 together with a canonical isomorphism
 \be  \xymatrix{
   {\u_\z} \ar[r]^<(0.3){can_\z} &   {\Delta_\z(M)} } \ee
 induced by the above choices of bases.\footnote{The choice of the basis $P_i^\vee$ of $T_{A^\vee}$ induces a map $\z^r\to T_{A^\vee}$ and, taking determinants, $can_{P^\vee}:\d_\z(\z^r)\to
 \d_\z(T_{A^\vee}).$    Similarly, we obtain canonical isomorphisms $can_{P^d},$ $can_{\gamma^+}$ and
 $can_{\delta}$ for $T^d_A,$ $T^+_B$ and $Lie_\z(A^\vee),$
 respectively. Set
 $can_\z:=can_{P^\vee}\cdot can_{P^d}^{-1}\cdot can_{\gamma^+}\cdot
 can_{\delta}^{-1}.$}

Define the period $\Omega^+_\infty(A)$ and the regulator $R_A$ of
$A$ to be the determinant of the maps $\alpha_M$ and $h$ with
respect to the  bases chosen above, respectively.  Then Conjecture
\ref{rationality} tells us that \be\label{canon}\zeta_\q(M)=\frac{
L^*(M)}{\Omega^+_\infty(A)\cdot R_A}\cdot can_\q,\ee where
$can_\q:\u_\q\to \Delta_\q(M)$ is induced from $can_\z$ by base
change. Indeed, we have  $\Omega^+_\infty(A)R_A=
(\vartheta_\infty)_\bbC\circ (can_\Q)_\bbC$ and
$L^*(M)=(\vartheta_\infty)_\bbC\circ (\zeta_\Q(M))_\bbC$ in
$\mathrm{Aut}_{\C_\bbC}(\u_\bbC)=\bbC^\times$ and thus
$\zeta_\Q(M)$ differs from $can_\Q$ by $\frac{
L^*(M)}{\Omega^+_\infty(A)\cdot R_A}.$

On the other hand,  using among others property h)  of the
determinant functor, Proposition \ref{resultsabelianvar} and the
identification $T_B^+\otimes_\z \zp\cong T_p^+$ (induced from
\eqref{B-l}) one easily verifies  that there is an isomorphism \be
\begin{split}
\Delta_\zp(M):=\Delta_\z(M)_\zp&\cong
\d_\zp(\r_f(\q,T_p))^{-1}\d_\zp(T_p^+)\d_\zp(Lie_\zp(A^\vee))^{-1}\\
&\quad\cdot
\d_\zp(\sha(A/\Q)(p))\d_\zp(A(\q)(p))^{-1}\d_\zp(A^\vee(\q)(p))^{-1},
\end{split}
\ee where $Lie_\zp(A^\vee):=Lie_\z(A^\vee)\otimes_\z \zp$
  is a $\zp$-lattice of
$t(M_p)\cong \H^1(A,\O_{A})\otimes_\q\qp.$ 

In order to compare this  with the integral structure
$\r_c(U,T_p)$ of $\r_c(U,M_p)$ we have to introduce the {\em local
Tamagawa numbers} $c_l(M_p)$ \cite[I \S 4]{fp}.

We first  assume $l\neq p.$ Then there is an exact sequence (cf.\
\cite[\S 4.2]{fp}) \be \xymatrix@C=0.5cm{
  0 \ar[r] & T^{I_l}_p \ar[rr]^{1-\phi_l} && T^{I_l}_p \ar[rr]^{ } && {\H^1_f(\ql,T_p) }\ar[rr]^{ } && {\H^1(I_l,T_p)_{tors}^{G_\ql} }\ar[r] & 0 }
  \ee
  which induces an isomorphism
  \begin{multline}
 \psi_l:\u_\zp\to \d_\zp([\xymatrix@C=0.5cm{
    T^{I_l}_p \ar[rr]^{1-\phi_l} && T^{I_l}_p }])\cong
  \d_\zp(\r_f(\ql,T_p))\d_\zp(\H^1(I_l,T_p)_{tors}^{G_\ql})^{-1}\\\cong\d_\zp(\r_f(\ql,T_p)).
  \end{multline}
  Here the first map arises as trivialization by the identity, the
  second comes from  the above exact sequence (interpreted as
  short exact sequence of complexes) while the last comes again from
  trivializing by the identity according to Remark
  \ref{finitemodules}.

  We define $c_l(M_p):=|\H^1(I_l,T_p)_{tors}^{G_\ql}|$ and remark
  that $(\psi_l)_\qp$
      differs from $\eta_l(M_p)$ precisely by the map $\xymatrix@C=0.5cm{
  {  \u_\qp}\ar[r]^<(0.1){acyc} & {\d_\zp(\d_\zp(\H^1(I_l,T_p)_{tors}^{G_\ql})_\qp} \ar@{=}[r]^<(0.4){def} & {\u_\qp,}
  }$
   which appealing to Remark  \ref{finitemodules} we also denote by $c_l(M_p).$ In other words we have \be(\psi_l)_\qp = c_l(M_p)\cdot\eta_l(M_p).\ee  Note
      also, that one has
  $|A^\vee(\ql)\otimes_\z\zp|=|P_l(M_p,1)|_p^{-1}\cdot c_l(M_p)$ and that
  $c_l(M_p)=1$ whenever $A$ has good reduction at $l.$

It can be shown \cite[Exp.\ IX,(11.3.8)]{SGA7} that $c_l$ is the
order of the $p$-primary part of the group of
$\mathbb{F}_l$-rational components $
(\mathcal{E}/\mathcal{E}^0)(\mathbb{F}_l)\cong\mathcal{E}(\mathbb{F}_l)/\mathcal{E}^0(\mathbb{F}_l)\cong
A(\qp)/A_0(\qp)$ of the special fibre $\mathcal{E}:=\A_{
\mathbb{F}_l}$ of the smooth (but not necessarily proper) Neron
model $\A$ of $A$ over $\z.$\footnote{The first isomorphism is a
consequence of the theorem of Lang \cite{lang} that the map
$x\mapsto \phi(x)x^{-1}$ on the $\overline{k}$-rational points of
a connected algebraic group over  a finite field $k$ (with
Frobenius $\phi$) is surjective.}


Now let $l=p.$ Similarly one defines maps (both depending on the
choice of $\delta$) \be \psi_p:\u_\zp\to
\d_\zp(\r_f(\qp,T_p))\d_\zp(Lie_\zp(A^\vee)),\ee \be
c_p(M_p):\u_\qp\to\u_\qp,\ee such that \be(\psi_p)_\qp =
c_p(M_p)\cdot\eta_p(M_p)\ee holds.\footnote{Assume that $A$ has
dimension $d$ and let $ \widehat{B}$ be the formal group of
$A^\vee$ over $\zp,$ i.e.\ the formal completion of $\B$ along the
zero-section in the fibre over $p.$ Note that $Lie_\zp(A^\vee)$
can be identified with the tangent space $t_{\widehat{B}}(\zp)$ of
$\widehat{B}$ with values in $\zp$ (a good reference for formal
groups is \cite{conrad-lieblich}). Furthermore we write
$\widehat{\mathbb{G}}_a$ for the formal additive group over $\zp,$
$\B^\circ$ and $ \widetilde{B}^\circ$ for the connected component
of the identity of $\B$ and its fibre $\widetilde{B}$ over $p,$
respectively, and $\Phi$ for the group of connected components of
$\widetilde{B}.$ Again by Lang's theorem we have
$\Phi(\fp)=\widetilde{B}(\fp)/ \widetilde{B}^\circ(\fp).$ Moreover
there are exact sequences
\[\xymatrix@C=0.5cm{
  0 \ar[r] & {\B^\circ(\zp)} \ar[rr]^{ } && {\B (\zp)} \ar[rr]^{ } && {\Phi(\fp)} \ar[r] & 0 }\]
and
\[\xymatrix@C=0.5cm{
 0 \ar[r] &  {\widehat{B}(\zp)} \ar[rr]^{ } && {\B^\circ(\zp)} \ar[rr]^{ } && {\widetilde{B}^\circ(\fp)} \ar[r] & 0. }\]
Now the logarithm map
\[\xymatrix@C=0.5cm{
  Lie_\zp(A^\vee)\supseteq (p\zp)^d=\widehat{\mathbb{G}}_a(\zp)^d  && {\widehat{B}(\zp)\subseteq A^\vee(\qp)^{\wedge p}=\H^1_f(\qp,T_p) }\ar[ll]^{log}_{\cong}
  }\]
induces the map $\psi_p$ by trivializing all finite subquotients
of the above line by the identity. Note that the first subquotient
on the left has order $p^d.$ Using \cite[ex.\ 3.11]{BK}, which
says that the Bloch-Kato exponential map coincides, up to the
identification induced by the Kummer map, with the usual
exponential map of the corresponding formal group,  it is easy to
see that $c_p(M_p):=\eta_p^{-1}\cdot
(\psi_p)_{\qp}=\bar{\eta_p}\circ (\psi_p)_{\qp}$ equals modulo
$\z^\times_p$
\begin{eqnarray*}
c_p(M_p)&=& p^{-d}|P(M_p,1)|_p^{-1}
\#\widetilde{B}^\circ(\fp)(p)\#\Phi(\fp)(p)=\#\Phi(\fp)(p),
\end{eqnarray*}
where we used the   relation
$|P(M_p,1)|_p=|P(M_l,1)|_p=p^{-d}\#\widetilde{B}^\circ(\fp)(p) .$
For elliptic curves this is well known \cite[appenndix \S
16]{silv}, the general case is an exercise using the description
of the reduction of abelian varieties in  \cite[Exp. IX]{SGA7}. }

Using the integral version of \eqref{triangle_c_f_f},  the maps
$\psi_l$ (analogously as $\eta_l$ for $\vartheta_p$
\eqref{vartheta})induce a canonical map \be
\kappa_p:\Delta_\zp(M)\to \d_\zp(\r_c(\q,T_p))\ee where all the
terms $\d_\zp(\sha(A/\Q)(p)),$ $ \d_\zp(A(\q)(p))^{-1}$ and $
\d_\zp(A^\vee(\q)(p))^{-1}$ are trivialized by the identity.
Hence, using again Remark \ref{finitemodules}, we have \be
\label{kappa} (\kappa_p)_\qp=
\frac{|\sha(A/\Q)|}{|A(\q)_{tors}||A^\vee(\q)_{tors}|}\prod
c_l(M_p) \cdot\vartheta_p\ee modulo $\z_p^\times.$ Since
$\zeta_\zp(T_p)$ equals $\kappa_p\circ can_\zp$ up to an element
in $\z^\times_p,$ it follows immediately from
\eqref{kappa},\eqref{canon} and \eqref{zetaintegral} that

\be \frac{ L^*(M)}{\Omega^+_\infty(A)\cdot R_A}\sim
\frac{|\sha(A/\Q)|}{|A(\q)_{tors}||A^\vee(\q)_{tors}|}\prod
c_l(M_p)\;\;\; \mbox{  mod } \z^\times_p. \ee

For all primes $p$ this implies the classical statement of the
BSD-Conjecture  up to sign (and a power of $2$ due to our
restriction $p\neq2$).

\section{The TNC - equivariant version}\label{ETNC}

The first equivariant version of the TNC with commutative
coefficients (other than number fields) was given by Kato
\cite{kato-LNM,kato1993} observing that classical Iwasawa theory
is, roughly speaking, nothing else than the ETNC for a "big"
coefficient ring. Inspired by Kato's work Burns and Flach
formulated an ETNC where the coefficients of the motive are
allowed to be (possibly non-commutative) finite-dimensional
$\Q$-algebras, using for the first time the general determinant
functor described in section \ref{det} and relative $K$-groups.
Their systematic approach recovers all previous versions of the
TNC and more over all central conjectures of Galois module theory.
 It were Huber and Kings \cite{hu-ki} who realized that the formulation of the
ETNC by relative $K$-groups is equivalent to the perhaps more
suggestive use of "generators," i.e.\ maps of the form $\u_R\to
\d_R(?)$ in the category $\C_R$ for various rings $R$ instead, see
also Flach's survey \cite[\S 6]{flach-survey}. They used this
approach to give - for motives of the form $M^*(1-k)$ with $k$ big
enough, i.e.\ with very negative weight - the first version of a
ETNC over general $p$-adic Lie extensions, which they call Iwasawa
Main Conjecture (while in this survey we reserve this name for
versions involving $p$-adic $L$-functions). While Burns and Flach
use "equivariant" motives and $L$-functions in their general
formalism, Fukaya and Kato realized that, at least for the
connection with Iwasawa theory which we have in mind, it is
sufficient to use non-commutative coefficients only for the Galois
cohomology, but to stick to number fields as coefficients for the
involved motives. In this survey we closely follow  their
approach.

To be more precise, consider for any motive $M$ the motive $
h^0(spec(F))\otimes M$ (both defined over $\Q$) for some finite
Galois extension $F$ of $\Q$ with Galois group $G=G(F/\Q).$ This
motive has a natural action by the group algebra $\z[G]$ and thus
will be of particular interest for Iwasawa theory where a whole
tower of finite extensions $F_n$ of $\Q$ is considered
simultaneously. Since there is an isomorphism of $K$-motives (for
$K$ sufficiently big)
\[h^0(spec(F))_K\otimes M\cong \bigoplus_{\rho \in \widehat{G}}  [\rho^*]^{n_\rho} \otimes M\]
where $\rho$ runs through all absolute irreducible representations
of $G$ and $n_\rho$ denotes the multiplicity with which it occurs
in the regular representation of $G$ on $K[G],$ it suffices - on
the complex side - to consider the collection of $K$-motives $
[\rho^*]  \otimes M$ and their $L$-functions or more precisely the
corresponding leading terms and vanishing orders. Indeed, the
$\bbC$-algebra $\bbC[G]$ can be identified with $\prod_{\rho \in
\widehat{G}} M_{n_\rho}(\bbC)$ and thus its first $K$-group
identifies with $\prod_{\rho \in \widehat{G}} \bbC^\times\cong
\mathrm{center}(\bbC[G]).$ In contrast, on the $p$-adic side, even
more when integrality is concerned, such a decomposition for
$\zp[G]$ is impossible in general.

This motivated Fukaya and Kato to choose the following form of the
ETNC. 
In fact, in order to keep the presentation  concise, we will only
describe a small extract of their complex and much more general
treatment.

Let $F$ be a $p$-adic Lie extension of $\Q$ with Galois group
$G=G(F/\Q).$ By $\La=\La(G)$ we  denote its Iwasawa algebra. For a
\Q-motive $M$ over \Q we fix a $G_\Q$-stable \zp-lattice $T_p$ of
$M_p$ and define a left $\La$-module
\[\T :=\La\otimes_\zp T_p\]
on which  $\La$ acts via multiplication on the left factor from
the left  while $G_\Q$ acts diagonally via $g(x\otimes y)=x
\bar{g}^{-1}\otimes g(y),$ where $\bar{g}$ denotes the image of
$g\in G_\Q$ in $G.$ This is a "big Galois representation" in the
sense of Nekovar \cite{nek}. Choose $S$ as in the previous section
and such that $\T$ is unramified outside $S$ and denote, for any
number field F', by $G_S(F')$ the Galois group of the maximal
outside $S$ unramified extension of $F'.$ Then by Shapiro's Lemma
 the cohomology of $\r(U,\T)$  for example is nothing else than
the perhaps more familiar $\La(G)$-module
$\H^i_{Iw}(F,T_p):=\varprojlim_{F'} \H^i(G_S(F'),T_p)$ where the
limit is taken with respect to corestriction and $F'$ runs over
all finite subextensions  of $F/\Q.$

Let $K$ be a finite extension of $\Q,$ $\lambda$ a finite place of
$K,$ $\O_\lambda$ the ring of integers of the completion
$K_\lambda$ of $K$ at $\lambda$ and assume that $\rho:G\to
GL_n(\O_{\lambda})$ is a continuous representation of $G$ which,
  for some suitable choice of a basis, is the $\lambda$-adic
realization $N_\lambda$ of a some $K$-motive $N.$ We also write
$\rho$ for the induced ring homomorphism $\La\to M_n(
\O_{\lambda})$ and we consider $\O_\lambda^n$ as a right
$\Lambda$-module via action by $\rho^t$ on the left,  viewing
$\O_\lambda^n$ as set of column vectors (contained in
$K_\lambda^n.)$ Note that, setting $M(\rho^*):= N^*\otimes M,$ we
obtain an isomorphism of Galois representations
\[\O_\lambda^n\otimes_\La \T\cong T_\lambda(M(\rho^*)),
\]
where $T_\lambda(M(\rho^*))$ is the $\O_\lambda$-lattice
$\rho^*\otimes T_p$ of $M(\rho^*)_\lambda$ and $\rho^*$ denotes
the contragredient(=dual) representation of $\rho.$


Now the equivariant version of conjecture \ref{integrality} reads
as follows

\begin{conj}[Equivariant Integrality; Fukaya/Kato]\label{equivintegrality}
There exists a (unique)\footnote{\label{uniqueness} In fact,
Fukaya and Kato assign such an isomorphism to each pair $(R,\T)$
where $R$ belongs to a certain class of rings containing the
Iwasawa algebras  for arbitrary $p$-adic Lie extensions of $\Q$ as
well as the valuation rings of finite extensions of $\qp$ and
where $\T$ is a projective $R$-module endowed  with a continuous
$G_\Q$-action. Then $\zeta_?(?)$ is supposed to behave well under
arbitrary change of rings for such pairs. Moreover they require
that the assignment $\T\mapsto \zeta_R(\T)$ is multiplicative for
short exact sequences. Only this full set of conditions   leads to
the uniqueness \cite[\S 2.3.5]{fukaya-kato}, while e.g.\ for
finite a group $G$ the map $K_1(\zp[G])\to K_1(\Q [G])$ need not
be injective and thus $\zeta_{\zp[G]}(?)$ might not be unique if
considered alone.} isomorphism
\[\zeta_\La(M):=\zeta_\La(\T): \u_\La \to \d_\La(\r_c(U,\T))^{-1}\]
with the following property:

For all $K, \lambda$ and $ \rho$ as above the  (generalized) base
change $\O_\lambda^n\otimes_\La -$ sends $\zeta_\La(M)$ to
$\zeta_{\O_\lambda}(T_\lambda(M(\rho^*))).$
\[\]

\end{conj}

Note that this conjecture assumes conjecture \ref{integrality} for
all $K$-motives $M(\rho^*)$ with varying $K.$ Furthermore, it is
independent of the choice of $S$ and of the lattices $T_p(M)$ and
$T_\lambda(M(\rho^*)).$

One obtains a slight modification - to which we will refer as the
{\em Artin-version }- of the above conjecture by restricting the
representations $\rho$ in question to the class of all Artin
representations of $G,$ (i.e.\ having finite image).  If $F/\Q$ is
finite, both versions coincide. Moreover, it is easy to
see\footnote{ If we assume that $K$ is big enough such that $K[G]$
decomposes completely into matrix algebras with coefficients in
$K,$ then the equivariant integrality statement (inducing
(absolute) integrality for $M(\rho^*)$ for all Artin
representations of $G$) amounts to an integrality statement for
the generator $\zeta_{K_\lambda[G]
}(M):=K_\lambda[G]\otimes_{\zp[G]}\zeta_{\zp[G]}(M)$   and thus
Burns and Flach's version for the $\Q$-algebra $K[G]$ with order
$\O_K[G].$ Using the functorialities of their construction
\cite[thm.\ 4.1]{bf}, it is immediate that taking norms leads to
the conjecture for the pair $(\Q[G],\z[G]).$ } that in this
situation the conjecture is equivalent to (the $p$-part of) Burns
and Flach's equivariant integrality conjecture \cite[conj.\ 6]{bf}
for the $\Q$-algebra $\Q[G]$ with $\z$-order $\z[G]$. Also,
$\T=\z_p[G]\otimes T_p(M)$ identifies with the induced
representation $\mathrm{Ind}_{G_\Q}^{G_F}T_p(M).$

Assume now that $F=\bigcup_n F_n$ is the union of finite
extensions $F_n$ of $\Q$ with Galois groups $G_n.$ Putting
$\zeta_{\Q_p[G_n]}(M)=\Q_p[G_n]\otimes_\Lambda \zeta_\Lambda(M)$
one recovers the "generator" $\delta_p(G_n,M,k)$\footnote{To be
precise, this is only morally true, since Huber and Kings take for
the definition of their generators the leading coefficients of the
modified $L$-function without the Euler factors in $S.$ It is not
clear to what extent this is compatible with our formulation
above} (for $k$ big enough) in \cite{hu-ki} as
$\zeta_{\Q_p[G_n]}(M^*(1-k)).$ Hence, up to shifting and Kummer
duality, the Artin-version of conjecture \ref{equivintegrality}
for $F$ is (morally) equivalent to \cite[Conj.\ 3.2.1]{hu-ki}.
Hence, using \cite[lem.\ 6.0.2]{hu-ki}  we obtain the following

\begin{prop}[Huber/Kings]
Assume Conjectures \ref{order}, \ref{rationality} and
\ref{integrality} for all $M(\rho^*)$ where $\rho$ varies over all
absolutely irreducible Artin representations of $G.$ Then the
existence  of $\zeta_{\La(G)}(M)$ satisfying the Artin-version of
Conjecture \ref{equivintegrality}  is equivalent to the existence
of $\zeta_{\La(G_n)}(M)$  for all $n.$
\end{prop}

In general, as remarked in footnote \ref{uniqueness},
$\zeta_{\zp[G_n]}(M)$ might not be unique (if it is considered
alone). But it is realistic to hope uniqueness for infinite $G$
(cf.\ \cite{kato-K1}) and then the previous zeta isomorphism would
be unique by the requirement
$\zeta_{\zp[G_n]}(M)=\zp[G_n]\otimes_{\La(G)}\zeta_{\La(G)}(M).$
Indeed, this is true at least if $G$ is big enough, see
\cite[prop.\ 2.3.7]{fukaya-kato}. Moreover, as Huber and Kings
\cite[\S 3.3]{hu-ki} pointed out, by twist invariance (over
trivializing extensions $F/\Q$ for a given motive) arbitrary
zeta-isomorphism $\zeta_\La(M)$ are reduced to those of the form
$\zeta_{\zp[G(F/\Q)]}(\Q)$  for the trivial motive and where $F$
runs through all finite extensions of $\Q.$

{\bf Question:} Does the Artin-version imply the full version of
Conjecture \ref{equivintegrality} ?

\section{The functional equation and
$\epsilon$-isomorphisms}\label{functional}

The $L$-function of a $\Q$-motive satisfies conjecturally a
functional equation, which we want to state in the following way
(to ease the notation we suppress the subscript $\Q$ in this
section)

\[L(M,s)=\epsilon(M,s) \frac{L_\infty(M^*(1),-s)}{L_\infty(M,s)} L(M^*(1),-s)\]

where the factor $L_\infty$ at infinity is built up by certain
$\Gamma$-factors and certain powers of $2$ and $\pi$ depending on
the Hodge structure of $M_B.$ The $\epsilon$-factor decomposes
into local factors \[\epsilon(M,s)=\prod_{v\in S}
\epsilon_v(M,s),\] where the definition for finite places is
recalled in footnotes \ref{epsilon_p} and \ref{epsilon_l};
$\epsilon_\infty(M,s)$ is a constant equal to a power of
$i.$\footnote{ We fix once and for all the complex period $2\pi
i,$ i.e.\ a square root of $-1,$ and, for every $l,$ the $l$-adic
period $t="2\pi i",$ i.e.\ a generator of $\z_l(1).$} We assume
this conjecture. Then, taking leading coefficients induces
\[L^*(M)=(-1)^{\eta}\epsilon(M)
\frac{L_\infty^*(M^*(1))}{L_\infty^*(M)} L^*(M^*(1))\] where
$\epsilon(M)=\prod\epsilon_v(M)$ with
$\epsilon_\nu(M)=\epsilon(M,0)$ and $\eta$ denotes the order of
vanishing at $s=0$ of the completed $L$-function
$L_\infty(M^*(1),s)L(M^*(1),s).$

\begin{ex}
For the motive $M=h^1(A)(1)$ of an abelian variety one has
$L_\infty(M,s)=L_\infty(M^*(1),s)=2(2\pi)^{-(s+1)}\Gamma(s+1),$
$L^*_\infty(M)=L^*_\infty(M^*(1))=\pi^{-1},$
$\epsilon_\infty(M)=-1$ and $\eta=0.$
\end{ex}

It is in no way obvious that the ETNC is compatible with the
functional equation and Artin-Verdier/Poitou-Tate duality. The
following discussion is a combination and reformulation of
\cite[section 5]{bf} and \cite[Appendix C]{perrin2000}. In order
to formulate the precise condition under which the compatibility
holds we first return to the absolute case and define "difference"
terms
\[L_{dif}^*(M):=L^*(M)L^*(M^*(1))^{-1}=(-1)^{\eta}\epsilon(M)
\frac{L_\infty^*(M^*(1))}{L_\infty^*(M)}
\]
and
\[\Delta_{dif}(M):=\d_\Q(M_B)\d_\Q(M_{dR})^{-1}.\]
We obtain an  isomorphism
\[\vartheta^{PD}:\Delta_\Q(M)\cdot\Delta_\Q(M^*(1))^*\cong  \Delta_{dif}(M) \] which arises from
the mutual cancellation of the terms arising from motivic
cohomology, the following isomorphism
\[M_B^+\oplus (M_B^*(1)^+)^*\cong M_B^+\oplus M_B(-1)^+\cong
M_B,\] where the last map is $(x,y)\mapsto x+ 2\pi i y,$ and from
the Poincare duality exact sequence
\be\label{PD}\xymatrix@C=0.5cm{
  0 \ar[r] & (t_{M^*(1)})^* \ar[rr]^{ } && M_{dR} \ar[rr]^{ } && t_M \ar[r] & 0
  .}\ee
On the other side define an isomorphism
\[\vartheta^{dif}_\infty: \Delta_{dif}(M)_{\mathbb{R}}\cong \u_{\mathbb{R}}\]
applying the determinant to \eqref{B+-dR} and to the following
isomorphism \be\label{infty+-}(\bbC\otimes_\Q M_B)^+=(\R\otimes_\Q
M_B^+)\oplus (\R(2\pi i)^{-1}\otimes_\Q M_B^-)\cong
(M_B)_{\mathbb{R}}\ee where the last map is induced by $\R(2\pi
i)^{-1}\to \R, x\mapsto 2\pi i x.$

Due to the autoduality of  the exact sequence of Conjecture
\ref{motinfty} (see \cite[lem.\ 12]{bf}) we have a commutative
diagram
\[\xymatrix{
  {\Delta_\Q(M)\cdot\Delta_\Q(M^*(1))^*} \ar[d]_{\vartheta_\infty(M)\cdot \overline{\vartheta_\infty(M^*(1))^*}} \ar[r]^<(0.3){\vartheta^{PD}_\mathbb{R}} &{ \Delta_{dif}(M)} \ar[d]^{\vartheta^{dif}_\infty} \\
  {\u_\mathbb{R} } \ar[r]^{\id_{\u_\mathbb{R}} } & {\u_\mathbb{R}}   }\]

Thus we obtain the following

\begin{prop}[Rationality]\label{difrationality} Assume that Conjecture
\ref{rationality} is valid for the $\Q$-motive $M.$ Then it is
also valid for its Kummer dual $M^*(1)$ if and only there exists a
(unique) isomorphism \[\zeta^{dif}(M):\u_\mathbb{Q}\to
\Delta_{dif}(M)\] such that we have
\[\xymatrix{
   L_{dif}^*(M):\u_\bbC\ar[rr]^{\zeta^{dif}(M)_\bbC} &   & {\Delta_{dif}(M)_\bbC} \ar[rr]^{(\vartheta^{dif}_\infty)_\bbC} &  & {\u_\bbC}   }\]
\end{prop}

Putting    $t_H(M):=\sum_{r\in \z} rh(r) $ with $h(r):=\dim_K
gr^r(M_{dR})(=\dim_{K_\lambda}gr^r (D_{dR}(M_\lambda)))$ for a
$K$-motive $M$ and noting that $t_H(M)=t_H(\det(M)),$ we have in
fact the following

\begin{thm}[Deligne {\cite[thm.\ 5.6]{deligne-L}}, Burns-Flach
{\cite[thm.\ 5.2]{bf}}]\label{detM} If the motive $\det(M)$ is of
the form $\Q(-t_H(M))$ twisted by a Dirichlet character, then
$\zeta^{dif}(M)$ exists.
\end{thm}

Deligne \cite{deligne-L} conjectured that the condition of the
theorem is satisfied for all motives. It is known to hold in all
examples A)-E).

See \eqref{ratdif} below for the rationality statement which is
hidden in the formulation of this theorem. Now we have to check
the compatibility with respect to the $p$-adic realizations. To
this aim we define the isomorphism
\[\vartheta_p^{dif}:\Delta_{dif}(M)_\qp\cong\d_\qp(M_p)\cdot\prod_{S\setminus S_\infty} \d_\qp(\r(\ql,M_p))\]
 as follows: Apply the determinant to \eqref{B-l} and multiply the
 resulting isomorphism with \[\id_{\d(M_{dR})^{-1}}\cdot \prod_{l\in S\setminus S_\infty}\Theta_l(M_p):\d(M_{dR})^{-1}\to \prod_{S\setminus S_\infty} \d_\qp(\r(\ql,M_p))\]
where $\Theta_l(M_p)=\eta_l(M)\cdot \eta_l(M^*(1))$ is defined in
the appendix \ref{appendix}.

On the other hand Artin-Verdier/Poitou-Tate Duality induces the
following isomorphism

\begin{gather*}
\d_\qp(\r_c(U,M_p))^{-1}\cong\d_\qp(\r(U,M_p))^{-1}    \d_\qp(\bigoplus_{v\in S} \r(\Q_v, M_p))\phantom{mmmmmmmmmmmmm}   \\
\begin{split}
 &\cong \d_\qp(  \r_c(U,M_p^*(1))^*)  \d_\qp((M_p^*(1)^+)^*)  \prod_{v\in S}   \d_\qp(\r(\Q_v,
 M_p))\\
 &\cong
\d_\qp(  \r_c(U,M_p^*(1))^*)  \d_\qp(M_p(-1)^+)
\d_\qp(M_p^+)\prod_{l\in S\setminus S_\infty} \d_\qp(\r(\ql,
M_p)).
 \end{split}
\end{gather*}

Using the identification \be \label{p-adic+-} M_p^+\oplus
M_p(-1)^+=M_p^+\oplus M_p^-(-1)\cong M_p,\ee where the last map is
induced by multiplication with the $p$-adic period $t="2\pi
i":M_p^-(-1)\to M_p^-,$ we obtain
\[ \vartheta^{AV}_p:\d_\qp(\r_c(U,M_p))^{-1}\cdot \d_\qp(  \r_c(U,M_p^*(1))^*)^{-1}\cong \d_\qp(M_p)\cdot\prod_{S\setminus S_\infty} \d_\qp(\r(\ql,M_p)). \]

Again one has to check the commutativity of the following diagram
(cf.\ \cite[lem.\ 12]{bf})
\[\xymatrix{
{\Delta_\Q(M)_\qp\cdot\Delta_\Q(M^*(1))_\qp^*} \ar[d]_{ \vartheta_p(M)\cdot\overline{\vartheta_p(M^*(1))^*}} \ar[r]^{(\vartheta^{PD})_\qp} & {\Delta^{dif}(M)_\qp} \ar[d]^{\vartheta_p^{dif}} \\
   {\d_\qp(\r_c(U,M_p))^{-1}\cdot \d_\qp(  \r_c(U,M_p^*(1))^*)^{-1}}\ar[r]^<(0.1){\vartheta^{AV}_p} & { \d_\qp(M_p)\cdot\prod_{S\setminus S_\infty} \d_\qp(\r(\ql,M_p)). }   }\]

Note that analogous maps exist and analogous properties hold also
if we replace $M_p$ by a Galois stable $\zp$-lattice $T_p$ or even
by the free $\La$-module $\T.$ Thus we obtain the following

\begin{prop}[Integrality] Assume that Conjecture
\ref{integrality} is valid for the $\Q$-motive $M.$ Then it is
also valid for its Kummer dual $M^*(1)$ if and only there exists a
(unique) isomorphism \[\zeta^{dif}_\zp(T_p):\u_\zp\to {
\d_\zp(T_p)\cdot\prod_{S\setminus S_\infty} \d_\zp(\r(\ql,T_p)). }
\] which induces via $\qp\otimes_\zp -$ the following map
\[\xymatrix{
    {\u_\qp\ar[rr]^{\zeta^{dif}(M)_\qp}} &   & {\Delta^{dif}(M)_\qp} \ar[rr]^<(0.2){\vartheta^{dif}_p} &  & { \d_\qp(M_p)\cdot\prod_{S\setminus S_\infty} \d_\qp(\r(\ql,M_p)). }  }\]
If this holds we have, using the above identifications, the
functional equation
\[\zeta_\zp(T_p)=(\overline{\zeta_\zp(T_p^*(1))^*})^{-1} \cdot \zeta^{dif}_\zp(T_p)\]
\end{prop}

Note that $(\overline{\zeta_\zp(T_p^*(1))^*})^{-1}$ is the same as
$ \zeta_\zp(T_p^*(1))^*\cdot \id_{\d_\zp(\r_{c}(U,T_p^*(1))^*)}$
according to Remark \ref{inverse}(i).  Needless to say that all
the above has an analogous version for $K$-motives whose
formulation we leave to the reader. Then it is clear how the
equivariant version of this proposition looks like:

\begin{prop}[Equivariant Integrality] Assume that Conjecture
\ref{equivintegrality} is valid for the $\Q$-motive $M.$ Then it
is also valid for its Kummer dual $M^*(1)$ if and only if there
exists a (unique) isomorphism \[\zeta^{dif}_\La(M):\u_\La\to {
\d_\La(\T)\cdot\prod_{S\setminus S_\infty} \d_\La(\r(\ql,\T)). }
\]

with the following property $(*):$

For all $K, \lambda$ and $ \rho$ as before Conjecture
\ref{equivintegrality} the (generalized) base change
$\O_\lambda^n\otimes_\La -$ sends $\zeta^{dif}_\La(M)$ to
$\zeta_{\O_\lambda}^{dif}(T_\lambda(M(\rho^*))),$ the analogue of
$\zeta^{dif}_\zp(T_p)$ for the $K$-motive  $M(\rho^*).$

If this holds we have  the functional equation
\[\zeta_\La(M)=(\overline{\zeta_\La(M^*(1))^*})^{-1} \cdot \zeta^{dif}_\La(M)\]
\end{prop}

Thus we formulate the
\begin{conj}[Local Equivariant Tamagawa Number Conjecture]\label{dif-conj}
The isomorphism $\zeta^{dif}_\La(M)$ in the previous Proposition
exists (uniquely).
\end{conj}

\subsection{$\epsilon$-isomorphisms}

One obtains a refinement of the above functional equation if one
looks more closely to which part of the Galois cohomology (and
comparison isomorphisms) the factors occurring in $L_{dif}^*(M)$
belong precisely. We first recall from \cite{perrin2000}
\be\label{perin}\frac{L_\infty^*(M^*(1))}{L_\infty^*(M)}=\pm
2^{d_{-}(M)-d_+(M)}(2\pi)^{-(d_-(M)+t_H(M))}\prod_{j \in
\z}\Gamma^*(-j)^{-h_j(M)}\ee where $\Gamma^*(-j)$ is defined to be
$\Gamma(j)=(j-1)!$ if $j>0$ and $\lim_{s\to j}
(s-j)\Gamma(s)=(-1)^{j}((-j)!)^{-1}$ otherwise.

The factor $(2\pi)^{-(d_-(M)+t_H(M))}$ arises as follows. Assume
for simplicity that $\det(M)=\Q(-t_H(M)).$ Then fixing a $\Q$-
basis $\gamma=(\gamma^+,\gamma^-)$ of $M_B$ and
$\omega=(\delta_M,\delta_{M^*(1)})$ of $M_{dR}$ which induce the
canonical basis (cf.\ example A)) of $\det(M)_B$ and
$\det(M)_{dR},$ respectively, gives rise to a map
\be\label{can}\xymatrix{
  {\u_\Q} \ar[r]^<(0.2){can_{\gamma,\omega}} & {\d_\Q(M_B)\d_\Q(M_{dR})^{-1}.}   }\ee
Base and change and the comparison isomorphism \eqref{B-dR} induce
\[\xymatrix{
  {\u_\bbC \ar[rr]^<(0.2){(can_{\gamma,\omega})_\bbC }} && {\d_\Q(M_B)_\bbC\d_\Q(M_{dR})^{-1}_\bbC \ar[rr]^<(0.3){\d(g_\infty)}} && {\u_\bbC}   }\]

whose value $\Omega_\infty^{dif}$ in $\bbC^\times$ is nothing else
than the inverse of the determinant over $\bbC$ of the comparison
isomorphism
\[\bbC\otimes_\Q \mathrm{det}(M)_B\to \bbC\otimes_\Q \mathrm{det}(M)_{dR}\]
 and thus $\Omega_\infty^{dif}=(2\pi i)^{-t_H(M)}.$
 But note that due to the definition of \eqref{infty+-} the above map  differs from
 $(\vartheta_\infty^{dif})_\bbC\circ (can_{\gamma,\omega})_\bbC $
 by the factor $(2\pi i)^{-d_-(M)}.$ Thus we obtain an explanation
 of the factor $(2\pi
 i)^{-(t_H(M)+d_-(M))}.$ Moreover, Theorem \ref{detM} and Proposition
\ref{difrationality} tell
 us that
 \be\label{ratdif}\frac{L^*_{dif}(M)}{(2\pi i)^{-d_-(M)}\Omega_\infty^{dif}}= \pm 2^{d_{-}(M)-d_+(M)}\frac{\epsilon_\infty(M)}{i^{-(t_H(M)+d_-(M))}}\prod_{l\in S\setminus S_\infty}\epsilon_l(M)\prod_{j \in
 \z}\Gamma^*(-j)^{-h_j(M)}\ee
is rational and that $\zeta_{dif}(M)$ is the map
$can_{\gamma,\omega}$ multiplied by this rational number.

The factor $2^{d_{-}(M)-d_+(M)}$ arises as quotient of the
Tamagawa factors of $M$ and $M^*(1)$ at infinity. One can either
cover it by defining $\Theta_\infty$ (see below) or changing  the
last map of the identification in \eqref{p-adic+-} as follows: on
the summand $M_p^+$ multiply with $2$ and on $M_p^-$ by
$\frac{1}{2}$ (as Fukaya and Kato do).

The  map \eqref{can} induces $\u_\qp\to \Delta_{dif}(M)_\qp\cong\d_\qp(M_p)\d_\qp(D_{dR}(M_p))^{-1}$ and furthermore
\[\xymatrix{
  {\u_{B_{dR}} \ar[rr]^<(0.3){(can_{\gamma,\omega})_{B_{dR}}}} & &{\Delta_{dif}(M)_{B_{dR}}} \ar[rr]^{\d_\qp(g_{dR})_{B_{dR}}} && {\u_{B_{dR}}}   }\]

whose value $\Omega_p^{dif}$ in $B_{dR}^\times$ is nothing else
than the inverse of the determinant over $B_{dR}$ of the comparison
isomorphism
\[B_{dR}\otimes_\qp D_{dR}(\mathrm{det}(M)_p)\to B_{dR}\otimes_\qp \mathrm{det}(M)_{p}\]
 and thus $\Omega_p^{dif}=(2\pi i)^{-t_H(M)},$ where we consider the $p$-adic period $ t=2\pi i$ as an
 element of $B_{dR}.$ Note that $  (B_{dR})^{I_p}=\qnr,$ the completion of the maximal
unramified extension $\mathbb{Q}_p^{nr}$ of $\qp.$ We need the
following

\begin{lem}[{\cite[prop.\
3.3.5]{fukaya-kato}},{\cite[C.2.8]{perrin2000}}]\label{qnr} The
map $\epsilon_p(M)\cdot\Omega_p^{dif}\cdot(can)_{B_{dR}}$ comes
from a map
\[\epsilon_{dR}(M_p):\u_{\qnr}\to\d_{\qnr}(M_p)\d_{\qnr}(D_{dR}(M_p))^{-1}.\]
Moreover, let $L$ be any finite extension of $\qp.$ Then a similar
statement holds for any finite dimensional $L$-vector space $V$
with continuous $G_\qp$-action instead of
$M_p\footnote{\label{epsilon_p}$\epsilon_p(V)=\epsilon(D_{pst}(V))$
where $D_{pst}(V)$ is endowed with the linearized action of the
Weil-group and thereby considered as a representation of the
Weil-Deligne group, see \cite[{\S}3.2]{fukaya-kato},   \cite{fp} or
\cite[appendix C]{perrin2000}. Furthermore, we suppress the
dependence of the choice of a Haar measure and of $t=2\pi i$ in
the notation. The choice of $t=(t_n)\in\z_p(1)$ determines a
homomorphism $\psi_p:\qp\to\bar{\qp}^\times$ with
$\ker(\psi_p)=\zp$ sending $\frac{1}{p^n}$ to
$t_n\in\mu_{p^n}.$}.$ We write $\epsilon_{dR}(V)$ for the
corresponding map, which is defined over $\widetilde{L}:=
\qnr\otimes_\qp L$ (see \cite[prop.\ 3.3.5]{fukaya-kato} for
details).
\end{lem}

For any $V$ as in the lemma we define an isomorphism
\[\epsilon_{p,L}(V):\u_{\widetilde{L}}\to
\big(\d_L(\r(\qp,V))\d_L(V)\big)_{\widetilde{L}}\]
 as product of $\Gamma_L(V):=\prod_\z \Gamma^*(j)^{-h(-j)},$ $\Theta_p(V)$ (see appendix \ref{Theta}) and $\epsilon_{dR}(V),$ where $h(j)=\dim_L gr^j D_{dR}(V).$

Now let $T$ be a Galois stable $\O:=\O_L$-lattice of $V$ and set $\widetilde{\O}:=W(\overline{\mathbb{F}_p})\otimes_\zp \O,$ where $W(\overline{\mathbb{F}_p})$ denotes the Wittring of $\overline{\mathbb{F}_p}.$ The following conjecture is a local integrality statement

\begin{conj}[Absolute $\epsilon$-isomorphism]\label{abs-eps-p}
There exists a (unique) isomorphism
\[\epsilon_{p,\O}(T): \u_{\widetilde{\O}}\to \big(\d_\O(\r(\qp,T))\d_\O(T)\big)_{\widetilde{\O}}\] with induces $\epsilon_{p,L}(V)$ by base change $L\otimes_\O -.$
\end{conj}

This conjecture, which is equivalent to conjecture $C_{EP}(V)$ in
\cite[III 4.5.4]{fp}, or more precisely its equivariant version
below  is closely related to the conjecture $\delta_\zp(V)$
\cite{perrin2000}   via the explicit reciprocity law
$R\acute{e}c(V),$ which was conjectured by Perrin-Riou and proven
independently by Benois \cite{benois}, Colmez \cite{colmez1998},
and Kurihara/Kato/Saito \cite{KKT}. In particular,  the above
conjecture is known for ordinary cristalline $p$-adic
representations \cite[1.28,C.2.10]{perrin2000} and for certain
semi-stable representations, see \cite{berger-tamagawa}.

To formulate an equivariant version, define
\[\widetilde{\La}:=\zpnr\kl G\kr=\projlim n \big(
W(\overline{\mathbb{F}_p})\otimes_\zp \zp[G/G_n]\big),\] where
$\zpnr=W(\overline{\mathbb{F}_p})$ denotes the ring of integers
 of $\qnr.$ We assume $L=\qp$ and set as before
$\T:=\La\otimes_\zp T$ (but later $T$ might differ from our global
$T_p$). We write $T(\rho^*)$ for the $\O$-lattice $\rho^*\otimes
T$ of $\rho^*\otimes V,$ which we assume de Rham.

\begin{conj}[Equivariant $\epsilon$-isomorphism]\label{eps-p}
There exists a (unique)\footnote{Again, Fukaya and Kato assign
such an isomorphism to each triple $(R,\T,t),$   where $R$ is as before,   $\T$ is a projective
$R$-module endowed  with a continuous $G_\qp$-action and $t$ is a generator of $\zp(1).$ Then
$\epsilon_{p,?}(?)$ is supposed to behave well under arbitrary change of
rings for such pairs. Moreover they require that the assignment
$\T\mapsto \epsilon_{p,R}(\T)$ is multiplicative for short exact
sequences, that it satisfies a duality relation when  replacing $\T$ by $\T^*(1),$ that the group $G_\qp^{ab}$ acts on a predetermined way (modifying $t$) compatible in a certain sense with the Frobenius ring homomorphism on $\widetilde{\La}$ induced from the absolute Frobenius of $\overline{\mathbb{F}_p}.$  Only this full set of conditions may lead to the uniqueness in general.} isomorphism
\[\epsilon_{p,\La}(\T): \u_{\widetilde{\La}}\to \big(\d_\La(\r(\qp,\T))\d_\La(\T)\big)_{\widetilde{\La}}\] such that for all $\rho:G\to GL_n(\O)\subseteq GL_n(L),$ $L$ a finite extension of $\qp$ with valuation ring $\O,$ we have
\[\O^n\otimes_\La \epsilon_{p,\La}(\T)= \epsilon_{p,\O}(T(\rho^*)).\]
\end{conj}

If $T=T_p\subseteq M_p$ is fixed we also write
$\epsilon_{p,\La}(M)$ for $\epsilon_{p,\La}(\T).$

Similarly we proceed in the case $l\neq p,$ formulating just one

\begin{conj}\label{eps-l}
There exists a (unique)\footnote{A similar comment as in the previous statement applies here.} isomorphism
\[\epsilon_{l,\La}(\T): \u_{\widetilde{\La}}\to \d_\La(\r(\ql,\T))_{\widetilde{\La}}\] such that for all $\rho:G\to GL_n(L),$ $L$ a finite extension of $\qp,$ we have
\[\O^n\otimes_\La \epsilon_{l,\La}(\T)= \epsilon_{l,\O}(T(\rho^*)).\]
Here $ \epsilon_{l,\O}(T(\rho^*)) $ is the analogue of the above
with respect to $\O$ instead of $\La$ and required to induce
\[L\otimes_\O \epsilon_{l,\O}(T(\rho^*)))=\Theta_l(V)\cdot
\epsilon_l(V(\rho^*))\] and its existence is part of the
conjecture\footnote{\label{epsilon_l} $\epsilon_l(V)=\epsilon( V)$
where $V$ is considered as representation of the Weil-Deligne
group of $\ql$ and where we suppress the dependence of the choice
of a Haar measure and of $t=2\pi i$ in the notation. The choice of
$t=(t_n)\in\z_l(1)$ determines a homomorphism
$\psi_l:\ql\to\bar{\ql}^\times$ with $Ker(\psi_l)=\z_l$ sending
$\frac{1}{l^n}$ to $t_n\in\mu_{l^n}.$ The formulation of this
conjecture is equivalent to \cite[conj.\ 3.5.2]{fukaya-kato} where
the constants $\epsilon_0$ are used instead of $\epsilon$ and
where $\theta_l$ does not occur. More precisely, our
$\epsilon_{l,\La}(\T)$ equals
$\epsilon_{0,\Lambda}(\ql,T,\xi)\cdot s_l(T)$ in
\cite[3.5.2,3.5.4]{fukaya-kato}}.
\end{conj}

For commutative $\La$ this conjecture was proved by S. Yasuda
\cite{yasuda} and it seems that he can extend his methods to cover
the non-commutative case, too.


If $T=T_p\subseteq M_p$ we also write $\epsilon_{l,\La}(M)$ for $\epsilon_{l,\La}(\T).$

Finally we set $\epsilon_{\infty,\La}(M)=\pm
2^{d_{-}(M)-d_+(M)}\frac{\epsilon_\infty(M)}{i^{-(t_H(M)+d_-(M))}}$
where the sign is  that which makes \eqref{perin} correct. The
following result is now immediate.

\begin{thm}[cf.\ {\cite[conj.\ 3.5.5]{fukaya-kato}}]
Assume Conjectures \ref{abs-eps-p}, \ref{eps-p} and \ref{eps-l}. Then Conjecture \ref{dif-conj}  holds, $\zeta^{dif}_\La(M)=\prod_{v\in S} \epsilon_{v,\La}(M)$ and we have the functional equation
\[\zeta_\La(M)=(\overline{\zeta_\La(M^*(1))^*})^{-1} \cdot \prod_{v\in S} \epsilon_{v,\La}(M)\]
\end{thm}

\section{$p$-adic $L$-functions and the Iwasawa main conjecture}
\label{p-adicL}

A $p$-adic $L$-function attached to a $\Q$-motive $M$ should be considered as a map on certain class of
representations of $G$ which interpolates the $L$-values of the  twists $M(\rho^*)$ at $0.$
The experience from those cases where such $p$-adic $L$-functions exist, shows that one has to
 modify the complex $L$-values by certain factors before one can hope to obtain a $p$-adic interpolation
 (cf.\ \cite{coates89,coates91} or \cite{perrin2000}
). The reason for this becomes clearer if one considers the Galois cohomology involved together with the
functional equation; in fact, that was the main motivation of the previous section.

In order to evaluate e.g.\ the $\zeta$-isomorphism or a
modification of it at a representation $\rho$ over a finite
extension $L$ of $\qp,$ one needs that the complex
$\rho\otimes_\La C,$   where $C$ is (a modification of)
$\r_c(U,\T),$ becomes acyclic: then the induced map $\u_L\to
\d_L(\rho\otimes_\La C)\to \u_L$   can be considered as value in
$K_1(L)=L^\times$ at $\rho.$

In general, $\r_c(U,\T)$ does not behave good enough and   will
have to be replaced  by some Selmer complex, which we well achieve
in two steps. This modification corresponds to a shifting of
certain Euler- and $\epsilon$-factors from one side of the
functional equation to the other such that  both sides are
balanced.

Though the following part of the theory holds in much greater
generality (e.g.\ in the {\em ordinary} good reduction case, but
not in the {\em supersingular} good reduction case) we just
discuss the case of abelian varieties in order to keep the
situation as concise as possible. Thus let $A$ be an abelian
variety over $\Q$ with good ordinary reduction at a fixed prime
$p\neq 2$  and set $M=h^1(A)(1)$ as before. Let $F_\infty$ be an
infinite $p$-adic Lie extension of $\Q$ with Galois group $G.$ For
simplicity we assume also that $G$ has no element of order $p,$
hence its Iwasawa algebra $\La=\La(G)$ is a regular ring.

Due to our assumption on the reduction type of $A,$ we have the
following fact: There is a unique $\qp$-subspace $\hat{V}$ of
$V=M_p$ which is stable under the action of $G_\qp$ and such that
\be\label{ordinary} D_{dR}(\hat{V})\cong D_{dR}(V)/D_{dR}^0(V).\ee
More precisely, $\hat{V}=V_p( \widehat{A^\vee})$ where
$\widehat{A^\vee}$ denotes the formal group  of the dual abelian
variety $A^\vee,$ i.e.\ the  formal completion of the Neron model
$\A/\zp$ of $A^\vee$ along the zero section of the special fibre $
\widetilde{\A}$. Then \eqref{ordinary} arises from the {\em unit
root splitting} $D_{dR}^0(V_p(\widetilde{\A}))\cong
D_{dR}^0(V)\subseteq D_{dR}(V)$ (see \cite[1.31]{nek-ht}) which is
induced from applying $D^0_{dR}(-)$ to the exact sequence of
$G_{\qp}$-modules
\[\xymatrix@C=0.5cm{
  0 \ar[r] & V_p( \widehat{A^\vee}) \ar[rr]^{ } && V_p(A^\vee) \ar[rr]^{ } && V_p(\widetilde{\A}) \ar[r] & 0 .}\]

Let $T$ be the $G_\Q$-stable $\zp$-lattice $T_p(A^\vee)$ of $V$
and set
\[ \hat{T}:=T\cap \hat{V}, \]
a  $G_\qp$-stable $\zp$-lattice of $\hat{V}.$ As before let $\T$
denote the big Galois representation $\La\otimes_\zp T$ and put
$\hat{\T}:=\La\otimes_\zp \hat{T}$ similarly. Then $\hat{\T}$ is
$G_\qp$-stable sub-\La-module of $\T.$ In fact, it is a direct
summand of $\T$ and  we have an isomorphism of $\La$-modules \be
\beta:\d_\La(\T^+)\cong \d_\La(\hat{\T}).\footnote{which arises as
follows: Choose a basis $\gamma^+=(\gamma_1^+,\ldots, \gamma_r^+)$
of $\H_1(A^\vee(\bbC),\Q)^+$ and $\gamma^-=(\gamma_1^-,\ldots,
\gamma_r^-)$ of $\H_1(A^\vee(\bbC),\Q)^-,$  which gives rise to a
$\zp$-basis of $T^+\cong (\H_1(A^\vee(\bbC),\z)\otimes_\z\zp)^+$
and $T^-\cong (\H_1(A^\vee(\bbC),\z)\otimes_\z\zp)^-$
respectively, where $r=d_+(M)=d_-(M).$
Then we obtain an isomorphisms \be \La^r \cong \T^+ \;\;\mbox{ and
}\;\; \phi:\d_\La(\La^r)\to \d_\La(T^+) \ee using the $\La$-basis
$\frac{1+\iota}{2}\otimes \gamma_j^+ + \frac{1-\iota}{2}\otimes
\gamma_j^-,$ $1\leq j\leq r,$ of $\T^+.$ On the other hand one can
choose   a $\Q$-basis $\delta=(\delta_1,\ldots,\delta_r)$   of
$Lie(A^\vee)$ (e.g.\ a $\z$-basis of $Lie_\z(A^\vee)$ as in
section \ref{BSD}) such that the isomorphism \be
\d_\qp(\Q_p^r)_{\qnr} \cong\d_\qp(\qp\otimes_\Q t_M)_{\qnr}\cong
\d_\qp(D_{dR}(\hat{V}))_{\qnr}\cong \d_\qp(\hat{V})_{\qnr},\ee
which is induced by $\delta,$ \eqref{dR-p}, \eqref{ordinary} and
$\epsilon_{dR}(\hat{V}))$ (according  to Lemma \ref{qnr}, but note
that $\epsilon_p(\hat{V})=1$ due to the good reduction), comes
from an isomorphism
\[\d_\zp(\mathbb{Z}_p^r)_{\zpnr}\cong \d_\zp(\hat{T})_{\zpnr},\]
where $\zpnr:= W(\fpbar).$ Then base change
$\widetilde{\La}\otimes_{\zpnr} -$ induces an isomorphism \be
\psi:\d_\La(\La^r)_{\widetilde{\La}}\cong
\d_\La(\hat{\T})_{\widetilde{\La}}.\ee Now
$\beta=\beta_{\gamma,\delta}$ is $\psi\circ\overline{\phi}.$}\ee

In this good ordinary case one can now first replace $\r_c(U,\T)$
by the Selmer complex $SC_U:=SC_U(\hat{\T},\T)$ (see \eqref{selmer_U}) which fits into the following distinguished triangle
\be\label{sc_u}\xymatrix{
  {\r_c(U,\T)}\ar[r]^{ } & {SC_U} \ar[r]^{ } & {\T^+\oplus\r(\qp,\hat{\T}) \ar[r]^{ }} &     }\ee
and thus induces an isomorphism
\[\d_\La (\r_c(U,\T))^{-1}\cong\d_\La (SC_U )^{-1}\d_\La (\T^+)\d_\La (\r(\qp,\hat{\T})).\]

The $p$-adic $L$-function $\L_U=\L_{U,\beta}(M,F_\infty/\Q)$
arises from the zeta-isomorphism $\zeta_\La(M)$ by a suitable
cancellation of the two last terms, compatible with the functional
equation. This is achieved by putting
\be\label{L}\L_U:=(\beta\cdot \id_{\d_\La(SC_U)^{-1}\d_\La(\hat{\T})^{-1}})\circ(\epsilon_{p,\La}(\hat{\T})^{-1}\cdot\zeta_\La(M)):\u_{\widetilde{\La}}\to \d_\La(SC_U)_{\widetilde{\La}}^{-1}.\ee

In order to arrive at a $p$-adic $L$-function which is independent of $U$ one has to replace $SC_U$ by another Selmer complex $SC:=SC(\hat{\T},\T) $ (\eqref{selmer}), which  fits into a distinguished triangle
\be \label{sc_u-sc}\xymatrix{
  SC_U \ar[r]^{ } & SC \ar[r]^{ } & {\bigoplus_{l\in S\setminus \{p,\infty\}} \r_f(\ql,\hat{\T})} \ar[r]^{ } &  ,  }\ee
  where for $l\neq p$ \be \label{r_flhat}\r_f(\ql,\hat{\T})\cong [\xymatrix{
    {\hat{\T}^{I_l}} \ar[r]^{1-\varphi_l} & {\hat{\T}^{I_l}}]  }\ee in the derived category.
    Replacing ${\hat{\T}^{I_l}}$ by a projective resolution if necessary and using the identity isomorphism  of it we obtain isomorphisms
\be \zeta_l(M)=\zeta_l(M,F_\infty/\Q): \u_\La\to
\d_\La(\r_f(\ql,\hat{\T}))^{-1}\ee and we define the $p$-adic
$L$-function \be\label{LU} \L=\L(M)=\L_U\cdot\prod_{l\in
S\setminus \{p,\infty\}} \zeta_l(M):\u_{\widetilde{\La}}\to
\d_\La(SC)_{\widetilde{\La}}^{-1}.\ee

  Let $\Upsilon$ be the set of all $l\neq p$ such that the ramification index of $l$ in $F_\infty/\Q$ is infinite. Note that $\Upsilon$ is empty if $G$ has a commutative open subgroup.

\begin{lem}\cite[prop.\ 4.2.14(3)]{fukaya-kato}\label{trivial}
$\hat{\T}^{I_l}=0$ and thus $\zeta_l(M)=1$ in $K_1(\La)$ for all $l$ in $\Upsilon.$
\end{lem}

Let us derive the {\em interpolation property} of $\L_U$ and $\L.$
Whenever $L^n\otimes_\La^\mathbb{L} SC_U$ is acyclic for a
continuous representation $\rho:G\to GL_n(\O_L),$ $L$ a finite
extension of $\qp,$ we obtain an element $\L_U(\rho)\in
\lnr^\times$ from the isomorphism
\[\xymatrix@C=0.5cm{
 { \u_{\tilde{L}}}\ar[rr]^<(0.2){L^n\otimes_\La \L_U} && {\d_{\tilde{L}} (L^n\otimes_\La^\mathbb{L}  SC_U)^{-1}}\ar[rr]^<(0.4){acyclic} && { \u_{\tilde{L}}}   }\]
 which via $K_1(\tilde{L})\to K_1(\lnr)$ can be considered as
 element of $\lnr^\times.$

Let $K$ be a finite extension of $\Q,$ $\rho:G\to GL_n(\O_{K})$ an Artin representation, $[\rho^*]$ the Artin motive corresponding to $\rho^*.$ Fix a place $\lambda$ of $K$ above $p,$   put $L:=K_\lambda$ and consider the $L$-linear representation of $G_\Q$ or its restriction to $G_{\ql}$
\[W:=M(\rho^*)_\lambda=[\rho^*]_\lambda\otimes_\qp M_p\]
and the $G_\qp$-representation
\[\hat{W}:=[\rho^*]_\lambda\otimes_\qp \hat{V}.\]

For a $G_\qp$-representation $V$ define
$P_{L,l}(V,u):=\det_L(1-\varphi_l u| V^{I_l})\in L[u]$ if $l\neq
p$ and $P_{L,p}(V,u):=\det_L(1-\varphi_p u| D_{cris}(V))\in L[u]$
otherwise.


Some conditions for acyclicity are summarized in the next

\begin{prop}
[{\cite[4.2.21, 4.1.6-8]{fukaya-kato}}] \label{acyclic} Assume the following conditions:

(i) $H^j_f(\Q,W)=H^j_f(\Q,W^*(1))=0$ for $j=0,1,$

(ii) $P_{L,l}(W,1)\neq 0$ for any $l\in\Upsilon$ (respectively for
any $l\in S\setminus \{p,\infty\}$).

(iii) $\{P_{L,p}(W,u)P_{L,p}(\hat{W},u)^{-1}\}_{u=1}\neq 0$ and $P_{L,p}(\hat{W}^*(1),1)\neq 0.$

Then the following complexes are acyclic:
$L^n\otimes_{\La,\rho}^\mathbb{L} SC$ (respectively $L^n\otimes_{\La,\rho}^\mathbb{L} SC_U$),  $\r_f(\ql, W)=L^n\otimes_{\La,\rho}^\mathbb{L}\r_f(\ql, \T),$ for any $l\in\Upsilon$ (respectively for any $l\in S\setminus \{p\}$). Furthermore, there is a quasi-isomorphism
\[\r(\qp,\hat{W})\to \r_f(\qp,W).\]
Finally, assuming Conjectures \ref{order} and \ref{motl}, $L_K(M(\rho^*),s)$ has neither zero or pole at $s=0.$
\end{prop}
 Henceforth we assume the conditions (i)-(iii).

We define $\Omega_\infty(M(\rho^*))\in \bbC^\times$  to be the
determinant of the period map
$\bbC\otimes_\mathbb{R}\alpha_{M(\rho^*)}$ with respect  to the
$K$-basis which arise from $\gamma$ (respectively $\delta$) and
the basis given by $\rho.$ It is easy to see that we have
\be\Omega_\infty(M(\rho^*))=\Omega_\infty^+(M)^{d_+(\rho)}\Omega_\infty^-(M)^{d_-(\rho)},\ee
where $d_\pm(\rho)=d_\pm([\rho])$ and $\Omega_\infty^\pm(M)$ is
the determinant of $\bbC\otimes_\Q M_B^\pm\cong \bbC\otimes_\Q
t_M$  with respect to the basis $\gamma^\pm$ and $\delta.$
Assuming Conjecture \ref{rationality} we have
\[\frac{L_K(M(\rho^*),0)}{\Omega_\infty(M(\rho^*))}\in K^\times.\]
We claim that, using Proposition  \ref{descent-prop}, the isomorphism $\L_U(\rho)$  
\begin{multline*}
\xymatrix@C=0.5cm{
  {\u }\ar[rrr]^<(0.3){\zeta_\La(M)(\rho)=}_<(0.3){\vartheta_\lambda\circ\zeta(M(\rho^*))} &&& {\d(\r_c(U,W))^{-1}} \ar[rrr]_<(0.2){ =\cdot\;\epsilon_{p,L}(\hat{W})^{-1}}^<(0.2){\cdot\;\epsilon_{p,\La}(\hat{\T})^{-1}(\rho)} &&& {\d(SC_U(\hat{W},W))^{-1}\d( \hat{W})^{-1}}\d( {W^+}) \ar[r] &}\\ \xymatrix@C=0.5cm{ {\ar[rr]^<(0.2){\id\cdot\beta(\rho)}}&& {\d(SC_U(\hat{W},W))^{-1}} \ar[rrr]^<(0.4){acyc}&& & {\u
  }}
\end{multline*}
 (we suppress for ease of notation the subscripts and remind the reader of our convention in Remark \ref{inverse}) is the product of the following automorphisms of $\u:$

 (1) $L_K(M(\rho^*),0) \Omega_\infty(M(\rho^*))^{-1},$

 (2) $\Gamma_L(\hat{W})^{-1}=\Gamma_\qp(\hat{V})^{-1},$

 (3) $\Omega_p(M(\rho^*))$ which is, by definition, the composite
 \be\xymatrix{
   {\d(\hat{W})}\ar[r]^{\cdot\;\epsilon_{dR}(\hat{W})^{-1}} &{\d(D_{dR}(\hat{W}))} \ar[r]^{\d(g_{dR}^t)}&{\d(t_{M(\rho^*)})}\ar[r]^<(0.2){\cdot\; can_{\gamma,\delta}}
   &{\d\big((M(\rho^*)_B^+)_L\big)} \ar[r]^<(0.3){\d(g_\lambda^+)}   & {\d(W^+)}\ar[r]^{\beta(\rho)}&  { \d(\hat{W})} }\ee
where we apply Remark \ref{inverse} to obtain an automorphism of
$\u,$\footnote{Using Remark \ref{inverse}(i) it is easy to see that this amounts to taking the
product of the following isomorphisms and identifying the target with $\u$ afterwards
\begin{align*}
&\xymatrix{{\u}\ar[rr]^<(0.2){can_{\gamma,\delta}} &   &
{\d\big((M(\rho^*)_B^+)_L\big)\d(t_{M(\rho^*)})^{-1},} } &
&\xymatrix{ {\u}\ar[rr]^<(0.25){\id_{-}\cdot{\d(g_\lambda^+)}  } & &{ \d(W^+) \d\big((M(\rho^*)_B^+)_L\big)^{-1},}}\\
&\xymatrix{ {\u}\ar[rr]^<(0.25){\id_{-}\cdot{\d(g_{dR}^t)}} &   &
{\d(D_{dR}(\hat{W}))^{-1}} \d(t_{M(\rho^*)}),} &
&\xymatrix{ {\u}\ar[rr]^<(0.3){\epsilon_{dR}(\hat{W})^{-1}} &   & {{ \d(\hat{W})^{-1}}\d(D_{dR}(\hat{W}))},   }\\
&\xymatrix{ {\u}\ar[rr]^<(0.3){\id_{-}\cdot\beta(\rho)} &   &{
\d(\hat{W})}{ \d(W^+)^{-1},} }
\end{align*}
where the identity maps are those of $\d\big((M(\rho^*)_B^+)_L\big)^{-1},$ $\d(D_{dR}(\hat{W}))^{-1}$ and $\d(W^+)^{-1},$ respectively.}

(4) $\prod_{S\setminus\{p,\infty\}} P_{L,l}(W,1):\xymatrix{
  {\u} \ar[rr]^<(0.3){\prod\eta_l(W)}   && {\prod \d(\r_f(\ql,W))} \ar[r]^<(0.4){acyc} & {\u }  }$ where the first map comes from the trivialization by the identity and the second from the acyclicity,

(5) $\{P_{L,p}(W,u)P_{L,p}(\hat{W},u)^{-1}\}_{u=1}:\xymatrix{
  {\u} \ar[r]^<(0.1){\eta_W\cdot\eta_{\hat{W}}^{-1}} &
  {\d(\r_f(\qp,W)) \d(\r(\qp,\hat{W}))^{-1}} \ar[r]^<(0.4){quasi} &
   {\u , }   }$ where we use that $t(W)=D_{dR}(\hat{W})=t(\hat{W})$ and the quasi-isomorphism mentioned in the above Proposition, and

(6) $P_{L,p}(\hat{W}^*(1),1):\xymatrix{
  {\u} \ar[rr]^<(0.2){\overline{(\eta_{\hat{W}^*(1)})^*}} && {\d(\r_f(\qp,\hat{W}^*(1)))} \ar[r]^<(0.4){acyc} & {\u,} }$ where we use that $t(\hat{W}^*(1))=D_{dR}^0(\hat{W})=0.$

For $\L$ we need beneath Lemma \ref{trivial} another

\begin{lem}[{\cite[lem.\ 4.2.23]{fukaya-kato}}]\label{Cf}
Let $l\neq p$ be not in $\Upsilon.$ Then
$L^n\otimes_{\La,\rho}^\mathbb{L} \r_f(\ql,\hat{\T})$ is acyclic
if and only if $P_{L,l}(W,1)\neq 0.$ If this holds then we have
$\zeta_l(M)(\rho)=P_{L,l}(W,1)^{-1}.$
\end{lem}

Thus we obtain the following

\begin{thm}[{\cite[thm.\ 4.2.26]{fukaya-kato}}]
Under the conditions (i)-(iii) from Proposition \ref{acyclic} and
assuming Conjecture  \ref{equivintegrality} for $M$ and Conjecture
\ref{eps-p} for $\hat{\T}$   the value $\L(\rho)$ (respectively
$\L_U(\rho)$) is
\begin{multline*}
 \frac{L_K(M(\rho^*),0)}{\Omega_\infty(M(\rho^*))}\cdot\Omega_p(M(\rho^*))\cdot\Gamma_\qp(\hat{V})^{-1} \cdot\\
 \cdot\{P_{L,p}(W,u)P_{L,p}(\hat{W},u)^{-1}\}_{u=1}\cdot
P_{L,p}(\hat{W}^*(1),1)\cdot \prod_B P_{L,l}(W,1),
\end{multline*}
 where
$B=\Upsilon\subseteq S\setminus\{p,\infty\}$ (respectively
$B=S\setminus \{p,\infty\}$).
\end{thm}

\begin{rem} Note that conditions (ii) and (iii) are satisfied in the case of an abelian variety with good ordinary reduction and $p.$ Furthermore, the quotient $\Omega_p(M(\rho^*))/\Omega_\infty(M(\rho^*)$ is independent of the
choice of basis $\gamma$ and $\delta.$ Also, it is easy to
see\footnote{Note that in the definition of $\beta$ and thus in
$\beta(\rho)$ the epsilon factor $\epsilon_p(\hat{V})$ in
$\epsilon_{dR}(\hat{V})$ equals $1$ and thus
$\Omega_p(M(\rho^*))=\beta(\rho)\circ
(\epsilon_p(\hat{W})^{-1}\cdot\overline{\beta(\rho)}).$} that for
some suitable choice we have
$\Omega_p(M(\rho^*))=\epsilon_p(\hat{W})^{-1}$ which, according to
standard properties of $\epsilon$-constants (cf.\
\cite[{\S}3.2]{fukaya-kato}) using that $\hat{V}$ is unramified as
module under the Weil-group, in turn is equal to
$\epsilon_p(\rho^*)^{-k} \cdot\nu^{-f_p(\rho)}$ where
$k=\dim_\Q(t_M)=\dim_\qp(\hat{V}),$
$\nu=\det_\qp(\varphi_p|D_{cris}(\hat{V}))$ and where $f_p(\rho)$
is the $p$-adic order of the Artin-conductor of $\rho.$ Due to the
compatibility conjecture $C_{WD}$ in \cite[2.4.3]{fontaine-ss},
which is known for abelian varieties (loc.cit., rem 2.4.6(ii)) and
for Artin motives, one obtains the $\epsilon$- and Eulerfactors
either from $D_{pst}(W)$ or from the corresponding $l$-adic
realisations with $l\neq p.$ Furthermore, we have
$P_{L,p}(W,1)\neq 0$ and $P_{L,p}(\hat{W},1)\neq 0$ for weight
reasons. Thus, noting that for abelian varieties
$\Gamma(\hat{V})=1,$ the above formula becomes

\begin{equation*}
 \tag{Int}\frac{L_{K,\Upsilon'}(M(\rho^*),0)}{\Omega_\infty(M(\rho^*))}\cdot \epsilon_p(\rho^*)^{-k} \cdot\nu^{-f_p(\rho)}
 \cdot \frac{P_{L,p}(\hat{W}^*(1),1)}{P_{L,p}(\hat{W},1) },
\end{equation*}

where $L_{K,\Upsilon'}$ denotes the modified $L$-function without
the Euler-factors in $\Upsilon':=\Upsilon\cup\{p\}.$\footnote{ In
order to compare this formula   with   (107) in \cite{cfksv} we
remark that, with the notation of (loc.\ cit.),
$u=\det(\phi_l|\widehat{V}(-1))= p\nu=p\omega^{-1}.$ Then by
\cite[rem.\ 4.2.27]{fukaya-kato} one has
$\epsilon(\rho^*)^{-d}\nu^{-f_p(\rho)}=\epsilon(\rho)^d
u^{-f_p(\rho)}$ (strictly speaking one has to replace the period
$t$ by $ -t$ in the second epsilon factor). !!!!! But it seems
that one has to interchange $\rho$ and $\widehat{\rho}$ on the
right hand side of (107)!!!!!!!!!!!!!!!!!! }
\end{rem}

\begin{proof} We consider the case $\L_U.$
First observe that due to the vanishing of the motivic cohomology the map
\[\zeta_K(M(\rho^*)):\u_K\to \d(\Delta(M(\rho^*)))=\d(M(\rho^*)_B^+)\d(t_{M(\rho^*)})^{-1}\] is just the map
$can_{\gamma,\delta}:\u\cong
\d(M(\rho^*)_B^+)\d(t_{M(\rho^*)})^{-1},$ induced by the bases
arising from $\gamma$ and $\delta,$    multiplied with
$L_K(M(\rho^*),0)  \Omega_\infty(M(\rho^*))^{-1}.$ Secondly, since
$\d(\r_f(\Q,W))=\u,$ the isomorphism
$\vartheta_\lambda(M(\rho^*)):\d(\Delta(M(\rho^*)))_L\cong
\d(\r_c(U,W))^{-1}$ corresponds up to  the identification
$\d(M(\rho^*)_B^+)_L\cong \d(W^+)$ to the product of
\[\xymatrix{
  {\d(t_{M(\rho^*)})^{-1}_L} \ar[rr]^{\overline{\d(g_{dR}^t)}^{\;-1}} && {\d(t(W))^{-1}} \ar[rr]^{\cdot\eta_\lambda(W)} && {\d(\r_f(\qp,W))}   }\] with
    $\prod_B P_{L,l}(W,1).$ Thirdly, the contribution from $ \epsilon_{p,L}(\hat{W})$ is  $\eta_p(\hat{W})\cdot \overline{\eta_p(\hat{W}^*(1))^*}\cdot \Gamma_L(\hat{W})\cdot \epsilon_{dR}(\hat{W})$ up to the canonical
    local duality isomorphism.  Together with $\beta(\rho)$ we thus obtain all the factors (1)-(6) above. To finish the proof in the case $\L$ use Lemmata \ref{trivial}, \ref{Cf} and \eqref{descent}.
    \end{proof}

\subsection{Interlude - Localized $K_1$}

The following construction of a localized $K_1$ is one of the
differences to the approach of Huber and Kings
\cite{hu-ki}\footnote{Instead of the localised $K_1$ they work
with $K_1$ of the ring $\projlim n \qp[G/G_n],$ which occurs in
the context of distributions, see \cite{colmez1998}.}. For a
moment let $\La$ be an arbitrary ring with unit and let $\Sigma$
be a full subcategory of $C^p(\La)$ satisfying (i) if $C$ is
quasi-isomorphic to an object in $\Sigma$ then it belongs to
$\Sigma,$ too, (ii) $\Sigma$ contains the trivial complex, (iii)
all translations of objects in $\Sigma$ belong again to $\Sigma$
and (iv) any extension $C$ in $C^p(\La)$ (by an exact sequence of
complexes) of $C',C''\in \Sigma$ is again in $\Sigma.$ Then Fukaya
and Kato construct a group $K_1(\La,\Sigma)$ whose objects are all
of the form $[C,a]$ with $C\in \Sigma$ and an isomorphism
$a:\u_\La\to\d_\La(C)$ (in particular, $[C]=0$ in $K_0(\La)$)
satisfying certain relations, see \cite[1.3]{fukaya-kato}. This
group fits into an exact sequence \be \xymatrix@C=0.5cm{
K_1(\La)\ar[rr]^{ } && K_1(\La,\Sigma) \ar[rr]^{\partial} &&
K_0(\Sigma) \ar[rr]^{ } && K_0(\La),
  }\ee
where $K_0(\Sigma)$ is the abelian group generated by $[[C]],$
$C\in \Sigma$ and satisfying certain relations. Here the first map
is given by sending the class of an automorphism $\La^r\to \La^r$
to $[[\La^r\to \La^r],can],$ where $can$ denotes the
trivialization of the complex $[\La^r\to \La^r]$ by the identity
according to Remark \ref{inverse}, $\partial$ maps $[C,a]$ to
$[[C]]$ while the last map is given by $[[C]]\mapsto [C].$  If $S$
is an left denominator set of $\La,$ $\La_S:= S^{-1}\La$ the
corresponding localization and $\Sigma_S$ the full subcategory of
$C^p(\La)$ consisting of all complexes $C$ such that $
\La_S\otimes_\La C$ is acyclic, then $K_1(\La,\Sigma_S)$ and
$K_0(\Sigma_S)$ can be identified with $K_1(\La_S)$ and
$K_0(S\mbox{-tor}^{pd}),$ respectively. Here $S\mbox{-tor}^{pd}$
denotes the category of $S$-torsion $\La$-modules with finite
projective dimension.

\subsection{ Iwasawa main conjecture I}

Let $\O$ be the ring of integers of the completion at any place
$\lambda$ above $p$ of the  maximal abelian outside $p$ unramified
extension  $F_\infty^{ab,p}$ of $\Q$ inside $F_\infty.$ Note that
the latter extension is finite because every non-finite abelian
$p$-adic Lie extension of $\Q$ contains the cyclotomic
$\zp$-extension, which is ramified at $p.$    Then by \cite[thm.\
4.2.26(2)]{fukaya-kato}  $\epsilon_{p,\La}(\hat{\T})$ and thus
$\L$ is already defined over $\La_\O:=\O\otimes_\zp \La(G)$
instead of $\widetilde{\La}.$

Now let $\Sigma=\Sigma_{SC}$ be the smallest full subcategory of
$C^p(\La_\O)$ containing $SC$ and satisfying the conditions
(i)-(iv) above. Then the evaluation  of $\L$ factorizes over its
class in $K_1(\La_\O,\Sigma)$ which we still denote by $\L.$ By
the construction of $\L$ we have

\begin{thm}[{\cite[thm.\ 4.2.22]{fukaya-kato}}]
Assume Conjectures \ref{equivintegrality} for $(M,\La)$ and
\ref{eps-p} for $(\hat{\T},\La).$ Then the following holds:

(i) $\partial(\L)=[[SC]]$

(ii) $\L$ satisfies the interpolation property (Int).
\end{thm}

{\bf Question:} If one knows the existence of $\L\in
K_1(\La_\O,\Sigma)$ with the above properties, what is missing to
obtain the zeta-isomorphism?


\subsection{Canonical Ore set}

Now assume that the cyclotomic $\zp$-extension $\Q_{cyc}$ is
contained in $F_\infty$ and set $H:=G(F_\infty/\Q_{cyc}).$ In this
situation  there exist a canonical left and right denominator set
of $\La_\O$
\[S^*=\bigcup_{i\geq 0} p^iS \] with \[S=\{\lambda\in
\La_\O|\La_\O/\La_\O\lambda \mbox{ is a finitely generated }
\La_\O(H)\mbox{-module}\}\] as was shown in \cite{cfksv}.

In this case we write $\M_H(G)$ for the category of $S^*$-torsion
modules and identify $K_0(\M_H(G))$ with $K_0(\Sigma_{S^*})$
recalling that $\La_\O$ is regular.

We write \[X= Sel(A/F_\infty)^\vee\] for the Pontryagin dual of
the classical Selmer group of $A$ over $F_\infty,$ see
\cite{cfksv}.

\begin{conj}[{\cite[conj.\ 5.1]{cfksv}}]
$X\in\M_H(G).$
\end{conj}

It is shown in \cite[prop.\ 4.3.7]{fukaya-kato} that the
conjecture is equivalent to   $SC$ belonging to $\Sigma_{S^*}.$ We
assume the conjecture. Observe that then
$\Sigma\subseteq\Sigma_{S^*},$ which induces a commutative diagram

\[ \xymatrix{
  K_1(\La_\O) \ar@{=}[d] \ar[r]^{ } & K_1(\La_\O,\Sigma) \ar[d]_{ } \ar[r]^{ \partial} & K_0(\Sigma) \ar[d]_{ } \ar[r]^{ } & K_0(\La_\O) \ar@{=}[d]^{} \\
  K_1(\La_\O) \ar[r]^{ } & K_1(({\La_\O})_{S^*}) \ar[r]^{\partial} & K_0(\M_H(G)) \ar[r]^{ } & K_0(\La_\O)
  }\]

{\bf Question:} Determine the (co)kernel of  the middle vertical
maps.

In \cite[\S 3]{cfksv}  it is explained how to evaluate elements of
$K_1(({\La_\O})_{S^*})$ at representations. By \cite[lem.
4.3.10]{fukaya-kato} this is compatible with the evaluation of
elements in $K_1(\La_\O,\Sigma).$ The following version of a Main
Conjecture was formulated in (loc.\ cit.).

\begin{conj}[Noncommutative Iwasawa Main Conjecture]\label{IMC}
There exist a (unique) element  $\L$ in $K_1(({\La_\O})_{S^*})$
such that

(i) $\partial \L=[X_\O]$ in $K_0(\M_H(G))$ and

(ii) $\L$ satisfies the interpolation property (Int).
\end{conj}

The connection with the previous version is given by the following

\begin{prop}
Let $F_\infty$ be e.g.\ $\Q(A(p))$ or
$\Q(\mu(p),\sqrt[p^\infty]{\alpha})$ for some $\alpha\in
\Q^\times\setminus\mu$ (false Tate curve)\footnote{ See
\cite[prop.\ 4.3.15-17]{fukaya-kato} for a more general
statement.}. Then \[[[X]]=[[SC]]\] in $K_0(\Sigma_{S^*}).$ In
particular, Conjecture \ref{IMC} is a consequence of Conjecture
\ref{equivintegrality} for $(M,\La)$ and Conjecture \ref{eps-p}
for $(\hat{\T},\La).$
\end{prop}

The advantage of the localisation $({\La_\O})_{S^*}$ relies on the
fact that one has an explicit description of its first $K$-group
since the natural map $({\La_\O})_{S^*}^\times\to
K_1(({\La_\O})_{S^*})$ induces quite often an isomorphism of the
maximal abelian quotient of $({\La_\O})_{S^*}^\times$ onto
$K_1(({\La_\O})_{S^*}),$ see \cite[thm.\ 4.4]{cfksv}. On the other hand,
the localized $K_1(\La_\O,\Sigma)$ exists without the assumption
that $G$ maps surjectively onto $\zp,$ e.g.\ if $G=SL_n(\zp).$
Also, if $G$ has $p$-torsion elements, i.e.\ if $\La_\O(G)$ is
{\em not} regular, one can still formulate  the Main Conjecture
using the complex $SC$ instead of the classical Selmer group $X$
(which could have infinite projective dimension).

{\bf Question:} To which extend does  a $p$-adic $L$-function $\L$
together with Conjecture \ref{IMC} determine the
$\zeta$-isomorphism  in Conjecture \ref{equivintegrality}? In
other words, does the Main conjecture imply the ETNC?

\section{Appendix: Galois cohomology}\label{appendix}

The main reference for this appendix is \cite[\S
1.6]{fukaya-kato}, but see also \cite{bf,bf3}. For simplicity we
assume $p\neq 2$ throughout this section. Let
$U=spec(\z[\frac{1}{S}])$ be a dense open subset of $spec(\z)$
where $S$ contains $S_p:=\{p\}$ and $S:=\{\infty\}$ (by abuse of
notation). We write $G_S$ for the Galois group of the maximal
outside $S$ unramified extension of $\Q.$ Let $X$ be a topological
abelian group with a continuous action of $G_S.$ Examples we have
in mind are
 $X=T_p, M_p, \T,$ etc.\ Using continuous cochains one defines a complex $\r(U,X)\footnote{For ease of notation we do
 not distinguish between complexes and their image in the derived category,
 though this is sometimes necessary in view of the correct use of the determinant functor and exact
 sequences of complexes.}$ whose cohomology is $\H^n(G_S,X).$ Then $\r_c(U,X)$ is defined by the  exact triangle
\be\xymatrix{
  {\r_c(U,X)}\ar[r]^{ } & {\r(U,X)}\ar[r]^{ } & {\bigoplus_{v\in S}\r(\Q_v,X)} \ar[r]^{ } &     }\ee
   where the $\r(\ql,X)$ and $\r(\R,X)$ denote the continuous cochain complexes calculating the local Galois groups $\H^n(\ql,X)$ and $\H^n(\R,X).$
   Its cohomology is concentrated in degrees $0,1,2,3.$

Let $L$ be a finite extension of $\qp$ with ring of integers $\O.$ Now we define the local and global {\em "finite parts"} for a
finite dimensional $L$-vector space $V$ with continuous
$G_{\Q_v}$- and $G_\Q$-action, respectively. For $\Q_v=\R$ we set
\[\mathrm{R\Gamma}_f(\bbR,V):=\mathrm{R\Gamma}(\bbR,V)\] while for a finite place
$\r_f(\ql,V)$ is defined as a certain subcomplex of $\r(\ql,V),$
concentrated in degree $0$ and $1,$  whose image in the derived
category is isomorphic to
\begin{gather}
 \label{rf-l}
  \r_f(\ql,V)\cong \begin{cases}
   [\xymatrix{ {V^{I_l}}\ar[r]^{1-\varphi_l} & {V^{I_l}} }]  & \text{ if $l\neq p,$ }\\
   [\xymatrix{
      D_{cris}(V) \ar[rr]^<(0.2){(1-\varphi_p,1)} && D_{cris}(V)\oplus D_{dR}(V)/D_{dR}^0(V)  }] & \text{ if $l=p.$}
  \end{cases}
    \end{gather}
       Here $\varphi_l$ denotes the geometric Frobenius (inverse of the arithmetic) and the induced map $D_{dR}(V)/D_{dR}^0(V)\to H^1_f(\qp,V)$ is called exponential map $exp_{BK}(V)$ of Bloch-Kato, where we write $\H^n_f(\ql,V)$ for the cohomology of $\r_f(\ql,V).$

Defining $\r_{/f}(\ql,V)$ as mapping cone \be
\label{r-/f}\xymatrix{
  {\r_f(\ql,V)} \ar[r]^{ } & {\r(\ql,V)} \ar[r]^{ } & {\r_{/f}(\ql,V)} \ar[r]^{ } &     }\ee
 we finally define $\r_f(\Q,V),$ whose cohomology is concentrated in degrees $0,1,2,3,$  as mapping fibre
  \be \xymatrix{
    {\r_f(\Q,V)} \ar[r]^{ } & {\r(U,V)} \ar[r]^{ } & {\bigoplus_{S\setminus S_\infty} \r_{/f}(\ql,V)} \ar[r]^{ } &    . }\ee
This is independent of the choice of $U.$ The octahedral axiom
induces an exact triangle \be \label{trian-c-f} \xymatrix{
  {\r_c(U,V)} \ar[r]^{ } & {\r_f(\Q,V)} \ar[r]^{ } & {\bigoplus_{S}\mathrm{R\Gamma}_f(\Q_v,V) }
  \ar[r]^{ } &  .   }\ee

\subsection{Duality}

Let $G,$ $\La=\La(G),$ $\T$ as in section \ref{ETNC}. By abuse of notation we write $-^*$ for both (derived) functors $\rhom_\La(-,\La)$ and $\rhom_{\La^\circ}(-,\La^\circ).$ Then Artin-Verdier/Poitou-Tate duality induces the existence of the following distinguished triangle in the derived category of $\La$-modules
\be \label{AVPT} \xymatrix@C=0.5cm{
   { \r_c(U,\T) }\ar[rr]^{ } && {\r(U,\T^*(1))^*[-3] }\ar[rr]^{ } && {\T^+} \ar[r] &   }\ee
   and similarly for $T$ (a Galois stable $\O$-lattice of $V$) and $V$ as coefficients (with   $\La$ replaced by $\O $ and  $L,$ respectively).\footnote{A more precise form to state the duality is the following. Let $\r_{(c)}(U,\T)$ be defined like $\r_{c}(U,\T)$  but using Tate cohomology $\widehat{\r}  (\R,\T)$ instead of the usual group cohomology $\r(\R,\T).$ Then one has isomorphisms \[ \r(U,\T^*(1))^*\cong\r_{(c)}(U,\T)[3]\cong \r(U,\T^\vee(1))^\vee\] where $-^\vee=\Hom_{cont}(-\qp/\zp)$ denotes the Pontryagin dual.}

For the finite parts one obtains from Artin-Verdier/Poitou-Tate  and local Tate-duality the following isomorphisms
\be \label{local-finite-dual}\r_f(\ql,V)\cong (\r(\ql,V^*(1))/\r_f(\ql,V^*(1)))^*[-2], \ee
\be \label{global-finite-dual}\r_f(\Q,V)\cong \r_f(\Q,V^*(1))^*[-3].\ee

Set $t(V):=D_{dR}(V)/D_{dR}^0(V)$ if $l=p$ and $t(V)=0$ otherwise.
Trivializing $V^{I_l}$ and $D_{cris}(V),$ respectively, in
\eqref{rf-l} by the identity induces, for each $l,$ an isomorphism
\begin{eqnarray}\label{eta}
\eta_l(V):& \u_{L} &\to
 \d_{L}(\mathrm{R\Gamma}_f(\Q_l,V))\d_L(t(V)).
\end{eqnarray}

Then, setting $D(V)=D_{dR}(V)$ if $l=p$ and $D(V)=0$ otherwise,
the isomorphism
\be\label{Theta}\Theta_l(V):\u_L\to\d_L(\r(\ql,V))\cdot
\d_L(D(V))\ee is by definition induced from $\eta_l(V)\cdot
\overline{(\eta_l(V^*(1))^*)}$
 followed by an isomorphism induced by   local duality \eqref{local-finite-dual} and using  the analogue $D_{dR}^0(V)=t(V^*(1))^*$ of \eqref{PD} if $l=p.$\footnote{  More explicitly, $\theta_p(V)$ is obtained from applying the determinant functor to the following exact sequence \[\begin{split}
&\xymatrix{
  0\ar[r]^{ } & {\H^0(\qp,V)} \ar[r]^{ } & D_{cris}(V) \ar[r]^{ } & D_{cris}(V)\oplus t(V) \ar[rr]^{exp_{BK}(V)} && {\H^1(\qp,V)} \ar[r]^{ } &     }\\
  &\xymatrix{
    \ar[rr]^<(0.3){exp_{BK}(V^*(1))^* } && D_{cris}(V^*(1))^*\oplus t(V^*(1))^*\ar[r]^{ } & D_{cris}(V^*(1))^* \ar[r]^{ } & {\H^2(\qp,V)} \ar[r]^{ } & 0  }
\end{split}\] which arises from joining the defining sequences of $exp_{BK}(V)$ with the dual sequence for $exp_{BK}(V*(1))$ by local duality \ref{local-finite-dual}. }

\subsection{Selmer complexes} For $l\neq p$ we define
$\r_f(\ql,\T)$ as in \eqref{rf-l} and $\r_{/f}(\ql,\T)$ as in
\eqref{r-/f} with $V$ replaced by $\T,$ see also \eqref{r_flhat}.
We do {\em not} define $\r_f(\qp,\T)$ since there is in general no
integral version of $D_{cris}(V).$

The Selmer complex $SC_U(\hat{\T},\T)$ is by definition the
 mapping fibre \be \label{selmer_U}\xymatrix@C=0.5cm{
    { SC_U(\hat{\T},\T)\ar[rr]^{ } }&& {\r(U,\T)} \ar[rr]^{ } && {\r(\qp,\T/\hat{\T})\oplus\bigoplus_{S\setminus (S_p\cup S_\infty)} \r(\ql,\T)} \ar[r] &  }\ee
while $SC(\hat{\T},\T)$ is the mapping fibre
\be \label{selmer}\xymatrix@C=0.5cm{
     {SC(\hat{\T},\T)}\ar[rr]^{ } && {\r(U,\T)} \ar[rr]^{ } && {\r(\qp,\T/\hat{\T})\oplus\bigoplus_{S\setminus (S_p\cup S_\infty)} \r_{/f}(\ql,\T)} \ar[r] &  }.\ee
Thus by the octahedral axiom one obtains the distinguished triangles \eqref{sc_u}, \eqref{sc_u-sc} and, using Artin-Verdier/Poitou-Tate duality,
\be \xymatrix@C=0.5cm{  {\r(U,\T^\vee(1))^\vee} \ar[rr]^{} && {SC(\hat{\T},\T)} \ar[rr]^{ } && {\r(\qp,\hat{\T})\oplus\bigoplus_{S\setminus (S_p\cup S_\infty)}\r_f(\ql,\T)} \ar[r] &   }. \ee

With the notation of section \ref{p-adicL} the Selmer complexes $ SC_U(\hat{W},W)$ and $SC(\hat{W},W)$ are defined analogously and satisfy analogous properties.

The following properties \cite[(4.2),propositions 1.6.5, 2.1.3 and
4.2.15]{fukaya-kato} are necessary conditions for the existence of
the zeta-iso\-morphism $\zeta_\La(\T)$ in Conjecture
\ref{equivintegrality} and the $p$-adic $L$-functions $\L_U$
\eqref{LU} and  $\L$ \eqref{L}.
\begin{prop} The complexes $\r_c(U,\T)$ and
$SC_U(\hat{\T},\T)$ are perfect\footnote{$\r_c(U,\T)$ is even
perfect for $p=2,$ this is the reason that it is better for the
formulation of the ETNC than $\r(U,\T).$} and in $K_0(\Lambda)$ we
have
\[ [\r_c(U,\T)]=[SC_U(\hat{\T},\T)]=0.\] If $G$ does not have
$p$-torsion, also $SC(\hat{\T},\T)$ is perfect and we have
$[SC(\hat{\T},\T)]=0.$
\end{prop}

\subsection{Descent properties}

 For the evaluation at representations one needs good descent properties of the complexes involved.

 \begin{prop}\label{descent-prop}\cite[prop.\ 1.6.5]{fukaya-kato}
 With the notation as in section \ref{p-adicL} we have canonical isomorphisms (for all $l$)
 \begin{align*}
L^n\otimes_{\La,\rho} \r(U,\T)&\cong \r(U,W), &  L^n\otimes_{\La,\rho} \r_{c}(U,\T)&\cong \r_{c}(U,W), \\
L^n\otimes_{\La,\rho} \r_{(c)}(U,\T)&\cong \r_{(c)}(U,W), & L^n\otimes_{\La,\rho} \r(\ql,\T)&\cong \r(\ql,W),\\
L^n\otimes_{\La,\rho}SC_U(\hat{\T},\T)&\cong SC_U(\hat{W},W).
 \end{align*}
 For $l\not\in\Upsilon\cup S_p$ we also have: $L^n\otimes_{\La,\rho} \r_f(\ql,\T)\cong
 \r_f(\ql,W).$
  \end{prop}

But note that the complex $\r_f(\ql,\hat{\T})$ for $l\in\Upsilon$
and thus $SC(\hat{\T},\T)$ does {\em not} descent like this in
general. Instead, according to \cite[prop.\ 4.2.17]{fukaya-kato}
one has a distinguished triangle \be
\label{descent}\xymatrix@C=0.5cm{
  L^n\otimes_{\La,\rho}SC(\hat{\T},\T)\ar[rr]^{ } && SC(\hat{W},W) \ar[rr]^{ } && {\bigoplus_{l\in \Upsilon}\r_f(\ql,W)} \ar[r] &   }.\ee

\INPUT{../bib/Xbib.bib}
\INPUT{etnc.bbl} 

\bibliographystyle{amsplain}
\bibliography{../bib/Xbib}
\end{document}